\PassOptionsToPackage{unicode}{hyperref}
\PassOptionsToPackage{hyphens}{url}
\PassOptionsToPackage{dvipsnames,svgnames,x11names}{xcolor}
\documentclass[
]{article}
\usepackage{xcolor}
\usepackage{amsmath,amssymb}
\usepackage[T1]{fontenc}
\usepackage[utf8]{inputenc}
\usepackage{textcomp} 

\usepackage{lmodern}
\makeatletter
\@ifundefined{KOMAClassName}{
  \IfFileExists{parskip.sty}{%
    \usepackage{parskip}
  }{
    \setlength{\parindent}{0pt}
    \setlength{\parskip}{6pt plus 2pt minus 1pt}}
}{
  \KOMAoptions{parskip=half}}
\makeatother
\makeatletter
\ifx\paragraph\undefined\else
  \let\oldparagraph\paragraph
  \renewcommand{\paragraph}{
    \@ifstar
      \xxxParagraphStar
      \xxxParagraphNoStar
  }
  \newcommand{\xxxParagraphStar}[1]{\oldparagraph*{#1}\mbox{}}
  \newcommand{\xxxParagraphNoStar}[1]{\oldparagraph{#1}\mbox{}}
\fi
\ifx\subparagraph\undefined\else
  \let\oldsubparagraph\subparagraph
  \renewcommand{\subparagraph}{
    \@ifstar
      \xxxSubParagraphStar
      \xxxSubParagraphNoStar
  }
  \newcommand{\xxxSubParagraphStar}[1]{\oldsubparagraph*{#1}\mbox{}}
  \newcommand{\xxxSubParagraphNoStar}[1]{\oldsubparagraph{#1}\mbox{}}
\fi
\makeatother

\usepackage{longtable,booktabs,array}
\usepackage{multirow}
\usepackage{calc} 
\usepackage{etoolbox}
\makeatletter
\patchcmd\longtable{\par}{\if@noskipsec\mbox{}\fi\par}{}{}
\makeatother
\IfFileExists{footnotehyper.sty}{\usepackage{footnotehyper}}{\usepackage{footnote}}
\makesavenoteenv{longtable}
\usepackage{graphicx}
\makeatletter
\newsavebox\pandoc@box
\newcommand*\pandocbounded[1]{
  \sbox\pandoc@box{#1}%
  \Gscale@div\@tempa{\textheight}{\dimexpr\ht\pandoc@box+\dp\pandoc@box\relax}%
  \Gscale@div\@tempb{\linewidth}{\wd\pandoc@box}%
  \ifdim\@tempb\p@<\@tempa\p@\let\@tempa\@tempb\fi
  \ifdim\@tempa\p@<\p@\scalebox{\@tempa}{\usebox\pandoc@box}%
  \else\usebox{\pandoc@box}%
  \fi%
}
\def\fps@figure{htbp}
\makeatother

\NewDocumentCommand\citeproctext{}{}

\makeatletter
 \let\@cite@ofmt\@firstofone
 \def\@biblabel#1{}
 \def\@cite#1#2{{#1\if@tempswa , #2\fi}}
\makeatother
\newlength{\cslhangindent}
\newlength{\csllabelwidth}
\newenvironment{CSLReferences}[2] 
 {\begin{list}{}{%
  \setlength{\itemindent}{0pt}
  \setlength{\leftmargin}{0pt}
  \setlength{\parsep}{0pt}
  \ifodd #1
   \setlength{\leftmargin}{\cslhangindent}
   \setlength{\itemindent}{-1\cslhangindent}
  \fi
  \setlength{\itemsep}{#2\baselineskip}}}
 {\end{list}}
\usepackage{calc}

\newcommand{\CSLLeftMargin}[1]{\parbox[t]{\csllabelwidth}{\strut#1\strut}}
\newcommand{\CSLRightInline}[1]{\parbox[t]{\linewidth - \csllabelwidth}{\strut#1\strut}}

\providecommand{\tightlist}{%
  \setlength{\itemsep}{0pt}\setlength{\parskip}{0pt}}

\usepackage{algorithm}
\usepackage{algpseudocode}
\usepackage{xcolor}
\usepackage{nicematrix}
\definecolor{modeone}{RGB}{102, 194, 165}
\definecolor{modetwo}{RGB}{252, 141, 98}
\definecolor{modethree}{RGB}{141, 160, 203}
\definecolor{modefour}{RGB}{231, 138, 195}
\usepackage{arxiv}
\usepackage{orcidlink}
\usepackage{amsmath}
\usepackage[T1]{fontenc}
\makeatletter
\@ifpackageloaded{caption}{}{\usepackage{caption}}
\AtBeginDocument{%
\ifdefined\contentsname
  \renewcommand*\contentsname{Table of contents}
\else
  \newcommand\contentsname{Table of contents}
\fi
\ifdefined\listfigurename
  \renewcommand*\listfigurename{List of Figures}
\else
  \newcommand\listfigurename{List of Figures}
\fi
\ifdefined\listtablename
  \renewcommand*\listtablename{List of Tables}
\else
  \newcommand\listtablename{List of Tables}
\fi
\ifdefined\figurename
  \renewcommand*\figurename{Figure}
\else
  \newcommand\figurename{Figure}
\fi
\ifdefined\tablename
  \renewcommand*\tablename{Table}
\else
  \newcommand\tablename{Table}
\fi
}
\@ifpackageloaded{float}{}{\usepackage{float}}
\floatstyle{ruled}
\@ifundefined{c@chapter}{\newfloat{codelisting}{h}{lop}}{\newfloat{codelisting}{h}{lop}[chapter]}
\floatname{codelisting}{Listing}

\makeatother
\makeatletter
\@ifpackageloaded{caption}{}{\usepackage{caption}}
\@ifpackageloaded{subcaption}{}{\usepackage{subcaption}}
\makeatother
\usepackage{bookmark}
\IfFileExists{xurl.sty}{\usepackage{xurl}}{} 
\hypersetup{
  pdftitle={Geometry Parameterisation via a Structural Modal Basis for Aerodynamic Shape Optimisation},
  pdfauthor={Alwin Wang; Robert Carrese; Luca Brown; Joel Collins; Pier Marzocca},
  pdfkeywords={Aircraft design, Shape optimisation, Modal
superposition, Geometry parameterisation, MDO},
  colorlinks=true,
  linkcolor={blue},
  filecolor={Maroon},
  citecolor={Blue},
  urlcolor={Blue},
  pdfcreator={LaTeX via pandoc}}

\newcommand{\runninghead}{A Preprint }
\renewcommand{\runninghead}{MPM Aerodynamic Shape Optimisation }
\title{Geometry Parameterisation via a Structural Modal Basis for
Aerodynamic Shape Optimisation}
\def\asep{\\\\\\ } 
\author{\textbf{Alwin Wang}\\\\Boeing Aerostructures Australia (BAA),
Melbourne, VIC 3207, Australia\\\\\asep\textbf{Robert Carrese}\\\\Boeing
Aerostructures Australia (BAA), Melbourne, VIC 3207,
Australia\\\\\asep\textbf{Luca Brown}\\\\Boeing Aerostructures Australia
(BAA), Melbourne, VIC 3207, Australia\\\\\asep\textbf{Joel
Collins}\\\\Boeing Aerostructures Australia (BAA), Melbourne, VIC 3207,
Australia\\\\\asep\textbf{Pier Marzocca}\\\\Sir Lawrence Wackett Defence
and Aerospace Centre, RMIT University, Melbourne, VIC 3083,
Australia\\\\}
\date{}
\begin{document}
\maketitle
\begin{abstract}
Aircraft outer mould lines are often heavily constrained early by
payload requirements, manufacturability limits, and low-observability
considerations. The remaining aerodynamic design space is typically
small, highly constrained, and difficult to represent using a concise
and physically meaningful set of design variables. A structural modal
parameterisation method (MPM) is introduced for such problems. The
method constructs a tunable pseudo-structure and solves its eigenvalue
problem to obtain a modal basis, from which a selected subset of modes
parameterises the geometry. Boundary conditions, stiffness distribution,
density, thickness, and added masses are treated as intentional
basis-design variables that shape the admissible deformation space prior
to optimisation. The resulting design coordinates are independent of any
specific flow solver, making the approach broadly applicable. For
tightly constrained aerodynamic shape optimisation problems, it provides
a practical balance of parameter conciseness, robust constraint
awareness, and geometric flexibility.
\end{abstract}
{\bfseries \emph Keywords}
\def\sep{\textbullet\ }
Aircraft design \sep Shape optimisation \sep Modal
superposition \sep Geometry parameterisation \sep 
MDO

\section{Nomenclature}\label{sec-nomenclature}

{\def\LTcaptype{none} 
\begin{longtable}[]{@{}
  >{\raggedright\arraybackslash}p{(\linewidth - 2\tabcolsep) * \real{0.1207}}
  >{\raggedright\arraybackslash}p{(\linewidth - 2\tabcolsep) * \real{0.8793}}@{}}
\toprule\noalign{}
\begin{minipage}[b]{\linewidth}\raggedright
Symbol
\end{minipage} & \begin{minipage}[b]{\linewidth}\raggedright
Meaning
\end{minipage} \\
\midrule\noalign{}
\endhead
\bottomrule\noalign{}
\endlastfoot
\(C\) & Pseudo-structure linear constraint operator (e.g.~fixed boundary
constraint) \\
\(E(x)\) & Pseudo-structure Young's modulus field for basis shape
tuning \\
\(\nu(x)\) & Pseudo-structure Poisson's ratio field for basis shape
tuning \\
\(\rho(x)\) & Pseudo-structure density field for basis shape tuning \\
\(K\) & Pseudo-structure stiffness matrix \\
\(\lambda_i\) & Pseudo-structure \(i\)th eigenvalue, with
\(\lambda_i = \omega_i^2\) \\
\(M\) & Pseudo-structure mass matrix \\
\(\omega_i\) & Pseudo-structure \(i\)th natural angular frequency \\
\(\phi_i\) & Pseudo-structure \(i\)th mode shape \\
\(\Phi_s\) & Matrix of selected modes,
\(\Phi_s = [\phi_{i_1}\ \phi_{i_2}\ \cdots\ \phi_{i_s}]\), expressed on
the optimisation geometry \\
\(P\) & Optional column-permutation matrix in the pivoted QR
factorisation \\
\(Q\) & Optional Euclidean-orthonormal basis from the pivoted QR
factorisation \\
\(R\) & Optional upper-triangular factor in the pivoted QR
factorisation \\
\(s\) & Number of selected modes \\
\(a\) & Vector of selected modal amplitudes \\
\(e\) & Prandtl span efficiency factor \\
\(x_0\) & Baseline optimisation geometry \\
\(x(a)\) & Updated geometry written in the selected modal basis \\
\end{longtable}
}

\section{Introduction}\label{sec-introduction}

Aircraft outer mould lines (OML) are often constrained early in the
development lifecycle by payload accommodation, manufacturability,
systems integration, certification requirements, and low-observability
(LO) {[}1,2{]}. As such, aerodynamic shape optimisation is rarely
performed over a large unconstrained design space. Instead, the
available geometric freedom is often limited to small, highly
constrained deformations, making it difficult to recover meaningful
aerodynamic improvements while satisfying all non-aerodynamic
requirements.

This creates a parameterisation problem before it creates an
optimisation problem. If too many geometric variables are exposed, the
search becomes expensive and the resulting shapes can become difficult
to constrain or interpret. If too few variables are used, the design
space may exclude the changes that matter aerodynamically {[}3,4{]}. For
tightly constrained problems, a useful basis therefore needs coordinates
that are concise enough to search, flexible enough to recover useful
global and local shapes, and constraint-aware enough that the resulting
motions remain admissible.

The structural modal parameterisation method (MPM) addresses that
requirement by constructing a tunable pseudo-structure, solving its
eigenvalue problem, and selecting a subset of modes. Boundary conditions
and coupling assumptions determine what motion is admissible before
optimisation begins, while stiffness, density, thickness, and added
masses prioritise admissible motions. The resulting coordinates are
smooth, coupled, and engineer-interpretable because they are derived
from natural deformation patterns rather than isolated geometric
handles. As a result, a large part of the admissibility logic is built
into the basis itself before any optimiser is called.

This basis occupies a practical niche between direct geometric controls,
deformation fields, and data-driven reduced bases. Compared with direct
or morphing-based parameterisations, protected regions, support-defined
motion, and regional coupling can often be imposed during basis
definition rather than managed later through large numbers of local
variables. Compared with reduced bases trained from sampled shapes, the
pseudo-structure can be re-tuned without first generating a shape
library. The coordinates remain solver-agnostic, and the practical
advantage is strongest when the required design changes are smooth,
tightly constrained, and awkward to express with conventional geometry
variables.

\section{Related Work and Motivation}\label{sec-related-work}

Parameterisation methods in aerodynamic optimisation are usually judged
by trade-off quality rather than by novelty alone. Sobester frames that
trade-off in terms of conciseness, robustness, and flexibility {[}3{]}.
The same tension is visible in two-dimensional airfoil benchmarking
{[}5{]} and becomes more severe once the problem shifts to tightly
constrained three-dimensional aircraft geometry {[}4,6{]}. The
literature is grouped here into four practical families: compact
analytic parameterisations, direct geometric parameterisations,
deformation and morphing methods, and reduced-order parameterisations.

Compact analytic parameterisations prescribe geometry directly through a
concise and usually interpretable variable set. Hicks-Henne bump
functions remain useful for local airfoil modification because they are
simple and differentiable {[}7{]}, while PARSEC and Class-Shape
Transformation (CST) move toward stronger geometric meaning and broader
smooth-shape coverage {[}8--10{]}. OpenVSP provides a compact parametric
representation for many three-dimensional planform configurations
{[}11{]}. These methods remain attractive when the required deformation
aligns with the assumed design variables, but they become restrictive
once the edit is strongly three-dimensional or several connected
surfaces must move together.

Direct geometric parameterisations begin from an existing geometry
surface description and manipulate that description directly. Bézier,
B-spline, and NURBS control structures are common because they are
smooth, familiar, and CAD-compatible, although variable count and
coupling can grow quickly on complicated three-dimensional surfaces
{[}3,5{]}. Feature-based CAD workflows preserve design intent, and
isogeometric analysis reduces the separation between geometry and
analysis descriptions {[}12,13{]}. However, once the aircraft design
constraints no longer align with the feature tree or control-net
structure, local aerodynamic correction and constraint management can
become costly.

Deformation and morphing methods modify an existing shape through a
separate motion field rather than the native geometry description.
Free-form deformation remains the standard example for smooth motion of
complex three-dimensional surfaces {[}14--18{]}, while radial basis
function interpolation is widely used when smooth surface and mesh
motion is required on arbitrary point sets {[}19--21{]}. Vertex Morphing
extends the same family toward very high nominal design freedom through
filtered surface-node motion {[}22,23{]}. These methods provide
substantial shape reach, but admissibility and geometric logic usually
still have to be imposed through separate constraints, filtering, or
restrictions.

Reduced-order parameterisations reduce dimensionality explicitly. Proper
orthogonal decomposition and singular value decomposition-based
geometric bases identify dominant directions from sampled shape
variation, while active-subspace methods depend more strongly on
sensitivities and response structure than on a large geometry library
{[}24--27{]}. Machine-learned latent representations extend the same
compression logic through trained low-dimensional shape spaces
{[}28--30{]}. These methods can be effective when representative data or
stable response structure exists, but for tightly constrained
proprietary aircraft components the required priors may be unavailable,
expensive to generate, or invalidated when the governing constraints
change.

\begin{figure}

\centering{

\includegraphics[width=0.7\linewidth,height=\textheight,keepaspectratio]{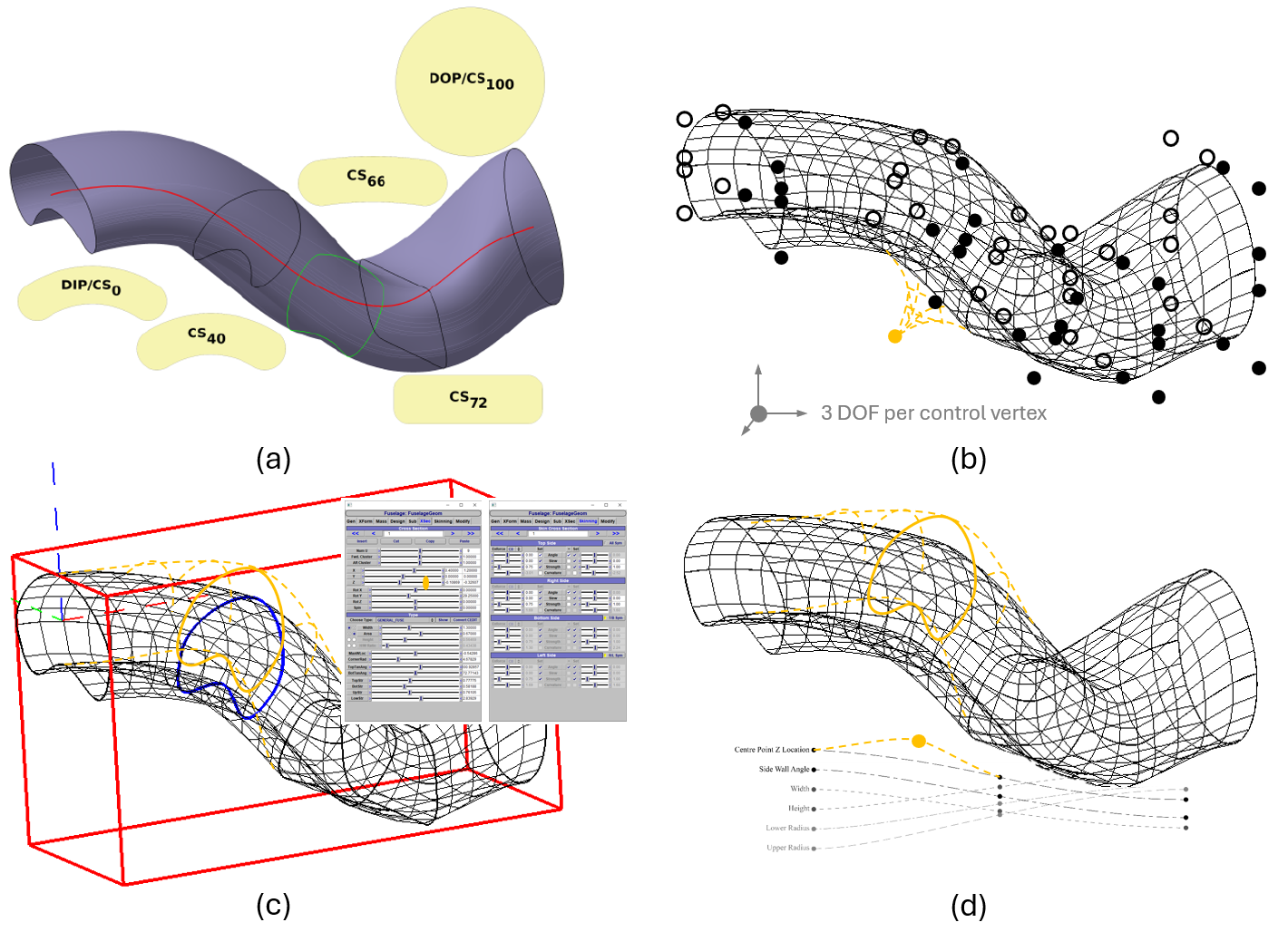}

}

\caption{\label{fig-traditional-parameterisation}Representative Military
Engine Intake Research Duct (MEIRD) geometry edits produced by
conventional parameterisation approaches. Orange overlays indicate
example surface deformations.}

\end{figure}%

\begin{figure}

\centering{

\includegraphics[width=0.55\linewidth,height=\textheight,keepaspectratio]{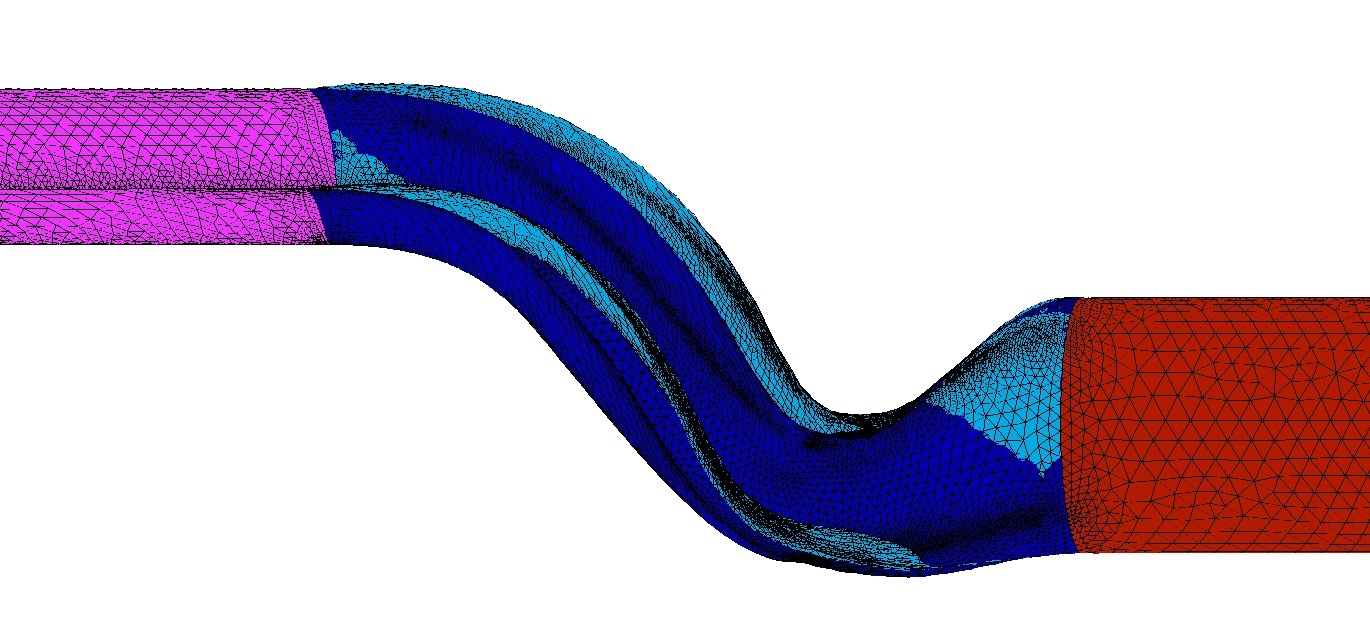}

}

\caption{\label{fig-cfd-mesh-morph}MEIRD CFD mesh morphing (light blue)
driven by selected pseudo-structural modes, illustrating transfer of
modal deformation to the aerodynamic mesh (dark blue).}

\end{figure}%

\begin{figure}

\centering{

\includegraphics[width=0.8\linewidth,height=\textheight,keepaspectratio]{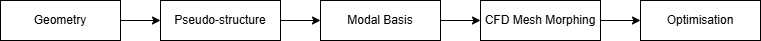}

}

\caption{\label{fig-flowchart}Practical workflow summary using MPM for
aerodynamic shape optimisation.}

\end{figure}%

Figure~\ref{fig-traditional-parameterisation} shows the MEIRD geometry
in panel (a) and examples of the design variables used by conventional
parameterisations: NURBS control points in panel (b), OpenVSP sections
in panel (c), and CAD guide and loft curves in panel (d).

Against that background, a pseudo-structure-derived basis occupies a
narrower and more practical niche. It seeks the variable economy of a
reduced-order representation without requiring a precomputed shape
library, while embedding a large part of the admissibility logic through
supports, coupling assumptions, and basis-shaping property edits.
Related work has demonstrated that structural motions can be transferred
to aerodynamic meshes for deformation and coupling {[}31{]}, as shown in
Figure~\ref{fig-cfd-mesh-morph} and Figure~\ref{fig-flowchart}. The step
taken here is to treat the tunable pseudo-structure itself as the
basis-design object, so that the modal deformation family can be revised
by editing the pseudo-structure before optimisation begins. That is
useful for tightly constrained aerodynamic geometry work because the
selected modes remain coupled, smooth, and engineer-interpretable while
still responding directly to packaging, attachment, and regional-motion
requirements.

\section{Modal Parameterisation
Method}\label{sec-modal-parameterisation-method}

\subsection{Define the Pseudo-Structure and Solve the Eigenvalue
Problem}\label{sec-define-pseudo-structure}

A family of admissible geometry changes is defined before any
optimisation objective is introduced. This is achieved by defining a
pseudo-structure over the geometry of interest and solving the
generalised eigenvalue problem

\begin{equation}\protect\phantomsection\label{eq-structural-eigenvalue-problem}{K \phi_i = \lambda_i M \phi_i}\end{equation}

where \(K\) is the pseudo-structural stiffness matrix, \(M\) is the
pseudo-structural mass matrix, \(\phi_i\) is the \(i\)th mode shape, and
\(\lambda_i = \omega_i^2\) is the associated eigenvalue.

The pseudo-structure is not required to represent the physical structure
of the final aircraft component. Rather, it is a basis-generation device
used to encode what motion should be forbidden, what motion should
remain admissible, and which natural deformation patterns are
prioritised among the candidate modes. In practice, this family of
admissible geometric deformations is defined through boundary
conditions, coupling, stiffness distribution, density distribution,
thickness distribution, and added masses. The pseudo-structure
discretisation also does not need to match the final structural-model
mesh. If built directly on the CFD surface mesh, the selected basis can
be applied to the associated CFD volume mesh through standard
mesh-morphing techniques {[}31{]}. These basis-design choices are then
fixed before the optimisation begins.

The imposed boundary conditions form an integral part of the
parameterisation definition because they determine what motion is
admissible before any aerodynamic objective is introduced. If \(C\)
denotes the linear constraint operator associated with fixed boundaries,
symmetry conditions, attachments, or coupling rules, then the admissible
pseudo-structural modes satisfy

\begin{equation}\protect\phantomsection\label{eq-structural-constraint-condition}{C \phi_i = 0}\end{equation}

This means that the selected basis already encodes those linear
geometric restrictions before optimisation begins. Structural-property
edits then steer which admissible motions appear early in the modal
ordering, but they do not replace the role of the constraints
themselves.

For distinct eigenvalues, the resulting mode shapes are orthogonal in
the mass metric. After mass normalisation, the candidate basis satisfies

\begin{equation}\protect\phantomsection\label{eq-structural-mass-stiffness-orthogonality}{\phi_i^\top M \phi_j =
\begin{cases}
0, & i\ne j \\
1, & i=j
\end{cases}
\quad
\text{where}
\quad
\phi_i^\top K \phi_j =
\begin{cases}
0, & i\ne j \\
\lambda_i, & i=j
\end{cases}}\end{equation}

This orthogonality is an important property, but it must be interpreted
in the appropriate space. The mode shapes are orthogonal in the
pseudo-structural mass metric \(M\). They should not be regarded as
inherently orthogonal under any metric that may be subsequently defined
on the geometry shape space. Additional conditioning steps may be
useful, but they are optional and are not taken here as the defining
feature of the proposed method. A selected subset of the mass-normalised
modes from one pseudo-structural eigenvalue problem is therefore full
rank and provides independent modal coordinates. Modes derived from
different pseudo-structures or from pseudo-structures with different
basis-design settings may nevertheless be combined to enrich the design
space. Such combined sets should be checked for linear dependence before
they are used as optimisation coordinates.

\subsection{Inspect, Tune, and Select the Subset of
Modes}\label{sec-tune-and-select-modes}

Once the candidate modes have been obtained, typically already
mass-normalised by commercial FEM solvers, they can be inspected
directly. The mode shapes show the admissible deformation behaviour
before any optimiser is called, so supports, stiffness, density,
thickness, and added masses can be adjusted until the candidate modal
family shows the desired global and local motion. This also supports
exploratory deformation studies before search begins, because the
candidate modal family can be inspected directly rather than inferred
from optimisation output.

\begin{equation}\protect\phantomsection\label{eq-modal-rayleigh-quotient}{\lambda_i = \frac{\phi_i^\top K \phi_i}{\phi_i^\top M \phi_i}}\end{equation}

The Rayleigh quotient in Equation~\ref{eq-modal-rayleigh-quotient} makes
the effect of basis-shaping edits easier to interpret. \(K\) acts as a
smoothing penalty for short-wavelength, high-curvature, or jagged
deformations. Reducing local stiffness therefore allows more localised
motion to appear earlier in the modal ordering. By contrast, \(M\)
provides the inertial weighting. Increasing local mass does not simply
force a node to move, but it can bias which regions participate more
strongly in the lower-order modes by changing the inertial cost of
moving that region.

Mode selection remains an engineering decision rather than a simple
consequence of eigenvalue order. Lower-order modes often carry useful
global behaviour, while higher-order modes often carry useful local
behaviour. The selected count \(s\) is therefore part of the basis
definition itself. After the selected modes have been expressed on the
optimisation geometry, they are assembled as

\begin{equation}\protect\phantomsection\label{eq-selected-structural-basis}{\Phi_s = \begin{bmatrix} \phi_{i_1} & \phi_{i_2} & \cdots & \phi_{i_s} \end{bmatrix}}\end{equation}

and the geometry update is written as

\begin{equation}\protect\phantomsection\label{eq-geometry-update}{x(a) = x_0 + \Phi_s a, \qquad a \in \mathbb{R}^s}\end{equation}

where \(a\) contains the selected modal amplitudes.

Algorithm \ref{alg:selected-basis-construction} summarises the
basis-definition loop. The pseudo-structure is revised until its modal
family meets the basis-design requirements. The selected modes are then
expressed on the optimisation geometry and their amplitudes used as
design variables.

\begin{algorithm}
\caption{Structural modal parameterisation using a tunable pseudo-structure}
\label{alg:selected-basis-construction}
\begin{algorithmic}[1]
\Require baseline geometry; pseudo-structure topology and discretisation; basis-design requirements
\State initialise the pseudo-structure definition
\Repeat
    \State assemble the pseudo-structural $K$ and $M$
    \State solve $K \phi_i = \lambda_i M \phi_i$
    \Comment mass-normalised eigenvectors
    \State inspect the candidate modes $\{\phi_i\}$
    \State update the pseudo-structure definition if required
    \Comment{constraints, stiffness, added masses}
\Until{the candidate modal family is satisfactory}
\State select $s$ useful modes, express them on the optimisation geometry, and assemble $\Phi_s$
\State use the amplitudes of $\Phi_s$ as design variables
\end{algorithmic}
\end{algorithm}

\subsection{Selected Basis and Optimisation
Coordinates}\label{sec-optimisation-coords}

The selected basis in Equation~\ref{eq-selected-structural-basis} is
used directly as the geometry parameterisation. If Euclidean-orthonormal
coordinates are preferred for a particular optimiser, a pivoted QR
factorisation may be applied to the selected basis,

\begin{equation}\protect\phantomsection\label{eq-selected-space-rotation}{\Phi_s P = Q R}\end{equation}

where \(P\) orders the selected modal columns before orthonormalisation,
\(Q\) has orthonormal columns, and \(R \in \mathbb{R}^{s \times s}\) is
upper triangular, so that

\begin{equation}\protect\phantomsection\label{eq-orthogonal-selected-basis}{Q^\top Q = I}\end{equation}

This factorisation preserves the selected deformation subspace because

\begin{equation}\protect\phantomsection\label{eq-rotated-selected-span}{\operatorname{span}(Q) = \operatorname{span}(\Phi_s)}\end{equation}

Pivoting changes the order of the selected modal columns, while QR
replaces the selected modal directions with Euclidean-orthonormal
coordinates. This optional transformation may reduce coordinate coupling
and improve numerical efficiency for some optimisation algorithms, but
its directions are generally less directly interpretable than the
selected pseudo-structural modes. It does not change the admissible
deformation family.

\section{MPM Applications}\label{sec-applications}

\subsection{Physics-derived Basis Functions: 2-Element Cantilever
Beam}\label{sec-cantilever-beam}

A simple 2-element Euler-Bernoulli cantilever beam problem was selected
as an analogy for the dihedral of a highly flexible transport aircraft
wing, such as the Boeing 787 at 1G cruise. The example is intentionally
small so that admissible and non-admissible dihedral shapes can be
inspected directly.

Two beam-shape libraries were generated from the 2-element beam model
with a clamped root and normalised tip displacement \(w_3 = 1\):

\begin{enumerate}
\def\labelenumi{\arabic{enumi}.}
\tightlist
\item
  Full arbitrary geometry family: all kinematically possible beam shapes
  including sagging, overshoot, and non-monotonically increasing shapes.
  These include shapes that would be unlikely as wing dihedral curves.
\item
  Desired dihedral geometry: non-negative displacement and non-negative
  slope, \(w(x) \ge 0\) and \(dw/dx \ge 0\), for \(0 \le x \le L\).
  These shapes would be representative candidates for wing dihedral
  curves.
\end{enumerate}

\begin{figure}

\centering{

\pandocbounded{\includegraphics[keepaspectratio]{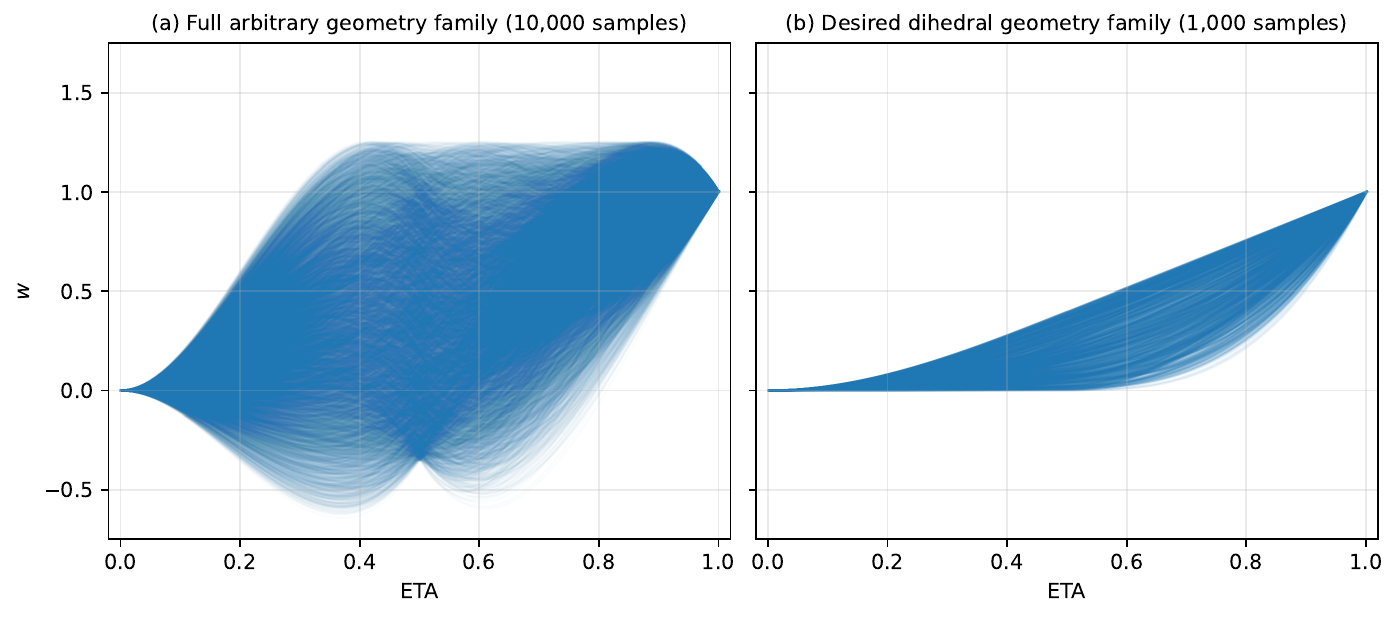}}

}

\caption{\label{fig-beam-geometry-families}Sampled shape libraries for
the 2-element cantilever beam with clamped root and normalised tip
displacement \(w_3 = 1\). Left: full arbitrary geometry family. Right:
desired dihedral family enforcing non-negative displacement and slope
over the span.}

\end{figure}

Both of these families were sampled uniformly for \(w_2\), \(\theta_2\),
and \(\theta_3\), and 10,000 curves were generated for the full
arbitrary family and 1,000 curves for the desired geometry family at 100
ETA (\(x/L\)) locations. These are shown in
Figure~\ref{fig-beam-geometry-families}.

For the pseudo-structural model, there are only four degrees of freedom
once the clamped root of the 2-element beam is enforced,

\begin{equation}\protect\phantomsection\label{eq-cantilever-reduced-coordinates}{q = \begin{bmatrix} w_2 & \theta_2 & w_3 & \theta_3 \end{bmatrix}^\top}\end{equation}

where \(w\) denotes transverse displacement and \(\theta\) denotes
section rotation. The standard Euler-Bernoulli element stiffness and
consistent mass matrices are

\begin{equation}\protect\phantomsection\label{eq-beam-element-stiffness-mass}{k^{(e)} = \frac{EI}{L_e^3}
\begin{bmatrix}
12 & 6L_e & -12 & 6L_e \\
6L_e & 4L_e^2 & -6L_e & 2L_e^2 \\
-12 & -6L_e & 12 & -6L_e \\
6L_e & 2L_e^2 & -6L_e & 4L_e^2
\end{bmatrix}
\quad
\text{and}
\quad
m^{(e)} = \frac{\rho A L_e}{420}
\begin{bmatrix}
156 & 22L_e & 54 & -13L_e \\
22L_e & 4L_e^2 & 13L_e & -3L_e^2 \\
54 & 13L_e & 156 & -22L_e \\
-13L_e & -3L_e^2 & -22L_e & 4L_e^2
\end{bmatrix}}\end{equation}

For illustrative purposes, \(A = L_e = E = I = 1\) and \(\rho = 42\)
were selected so that the resulting matrices and modes remained easy to
inspect. The dense assembled matrices show that the modal shapes couple
the \(w\) and \(\theta\) by construction rather than moving them
independently. The mass-normalised modal matrix \(\Phi\) satisfies the
expected mass orthogonality relation from
Section~\ref{sec-define-pseudo-structure}, although it is not
orthonormal itself.

\protect\phantomsection\label{cantilever-beam-modal-metrics}
\begin{equation}\protect\phantomsection\label{eq-cantilever-mass-orthogonality}{
\displaystyle
\underbrace{\begin{bNiceMatrix}[margin]
    \CodeBefore
    \rowcolor{modeone}{1}
    \rowcolor{modetwo}{2}
    \rowcolor{modethree}{3}
    \rowcolor{modefour}{4}
    \Body
    0.07 & 0.13 & 0.22 & 0.15 \\
    -0.16 & 0.05 & 0.22 & 0.53 \\
    0.02 & -0.94 & 0.25 & 1.18 \\
    0.1 & 1.07 & 0.41 & 3.98 \\
\end{bNiceMatrix}}_{\Phi^{\top}}
\underbrace{\begin{bmatrix}31.2 & 0.0 & 5.4 & -1.3\\0.0 & 0.8 & 1.3 & -0.3\\5.4 & 1.3 & 15.6 & -2.2\\-1.3 & -0.3 & -2.2 & 0.4\end{bmatrix}}_{M}
\underbrace{\begin{bNiceMatrix}[margin]
    \CodeBefore
    \columncolor{modeone}{1}
    \columncolor{modetwo}{2}
    \columncolor{modethree}{3}
    \columncolor{modefour}{4}
    \Body
    0.07 & -0.16 & 0.02 & 0.1 \\
    0.13 & 0.05 & -0.94 & 1.07 \\
    0.22 & 0.22 & 0.25 & 0.41 \\
    0.15 & 0.53 & 1.18 & 3.98 \\
\end{bNiceMatrix}}_{\Phi}
=
\underbrace{\begin{bmatrix}1 & 0 & 0 & 0\\0 & 1 & 0 & 0\\0 & 0 & 1 & 0\\0 & 0 & 0 & 1\end{bmatrix}}_{I}
}\end{equation}

For the following four-mode beam case, pivoted QR gives \(\Phi P = QR\),
where \(P\) orders the modal columns before orthonormalisation. In
\(Q = [q_1\ q_2\ q_3\ q_4]\), the subscript identifies the left-to-right
QR-column index after pivoting rather than an original modal index. The
subsequent comparison applies the same factorisation separately to each
greedily selected subset \(\Phi_s\).

\protect\phantomsection\label{cantilever-beam-modal-qr}
\begin{equation}\protect\phantomsection\label{eq-cantilever-pivoted-qr-factorisation}{
\displaystyle
\underbrace{\begin{bNiceMatrix}[margin]
    \CodeBefore
    \columncolor{modeone}{1}
    \columncolor{modetwo}{2}
    \columncolor{modethree}{3}
    \columncolor{modefour}{4}
    \Body
    0.1 & 0.02 & -0.16 & 0.07 \\
    1.07 & -0.94 & 0.05 & 0.13 \\
    0.41 & 0.25 & 0.22 & 0.22 \\
    3.98 & 1.18 & 0.53 & 0.15 \\
\end{bNiceMatrix}}_{\Phi P}
=
\underbrace{\begin{bNiceMatrix}[margin]
    \CodeBefore
    \columncolor{modefour}{1}
    \columncolor{modethree}{2}
    \columncolor{modetwo}{3}
    \columncolor{modeone}{4}
    \Body
    -0.03 & 0.0 & 0.75 & -0.67 \\
    -0.26 & -0.96 & -0.07 & -0.07 \\
    -0.1 & 0.13 & -0.66 & -0.73 \\
    -0.96 & 0.25 & 0.07 & 0.11 \\
\end{bNiceMatrix}}_{Q}
\underbrace{\begin{bmatrix}-4.14 & -0.92 & -0.54 & -0.2\\0.0 & 1.22 & 0.11 & -0.06\\0.0 & 0.0 & -0.23 & -0.09\\0.0 & 0.0 & 0.0 & -0.2\end{bmatrix}}_{R}
}\end{equation}

\begin{equation}\protect\phantomsection\label{eq-cantilever-pivoted-qr-orthogonality}{
\displaystyle
\underbrace{\begin{bNiceMatrix}[margin]
    \CodeBefore
    \rowcolor{modefour}{1}
    \rowcolor{modethree}{2}
    \rowcolor{modetwo}{3}
    \rowcolor{modeone}{4}
    \Body
    -0.03 & -0.26 & -0.1 & -0.96 \\
    0.0 & -0.96 & 0.13 & 0.25 \\
    0.75 & -0.07 & -0.66 & 0.07 \\
    -0.67 & -0.07 & -0.73 & 0.11 \\
\end{bNiceMatrix}}_{Q^{\top}}
\underbrace{\begin{bNiceMatrix}[margin]
    \CodeBefore
    \columncolor{modefour}{1}
    \columncolor{modethree}{2}
    \columncolor{modetwo}{3}
    \columncolor{modeone}{4}
    \Body
    -0.03 & 0.0 & 0.75 & -0.67 \\
    -0.26 & -0.96 & -0.07 & -0.07 \\
    -0.1 & 0.13 & -0.66 & -0.73 \\
    -0.96 & 0.25 & 0.07 & 0.11 \\
\end{bNiceMatrix}}_{Q}
=
\underbrace{\begin{bmatrix}1 & 0 & 0 & 0\\0 & 1 & 0 & 0\\0 & 0 & 1 & 0\\0 & 0 & 0 & 1\end{bmatrix}}_{I}
}\end{equation}

The basis vectors generated by this pseudo-structural model were
compared with those from the methods described in
Section~\ref{sec-related-work} for both geometry families shown in
Figure~\ref{fig-beam-geometry-families}:

\begin{itemize}
\tightlist
\item
  Direct beam-DOF parameterisation: using \(w_2\), \(\theta_2\),
  \(w_3\), and \(\theta_3\) directly
\item
  Raw MPM: selecting \(s\) raw structural modes \(\Phi_s\) from the
  mass-normalised modal matrix \(\Phi\) derived in
  Equation~\ref{eq-cantilever-mass-orthogonality}
\item
  Rotated MPM: selecting \(\Phi_s\) by the same greedy forward-selection
  rule at each selected count \(s\), then calculating \(Q^{(s)}\) by
  pivoted QR of that subset
\item
  Data-driven SVD: SVD was performed on the library of deflected shapes
  at the sampled ETA locations (a \(100 \times 10{,}000\) matrix for the
  full library and a \(100 \times 1{,}000\) matrix for the desired
  library)
\end{itemize}

Each \(Q^{(s)}\) was calculated independently from its selected
\(\Phi_s\) subset. A single QR factorisation of all four modes was not
used.

\begin{figure}

\centering{

\pandocbounded{\includegraphics[keepaspectratio]{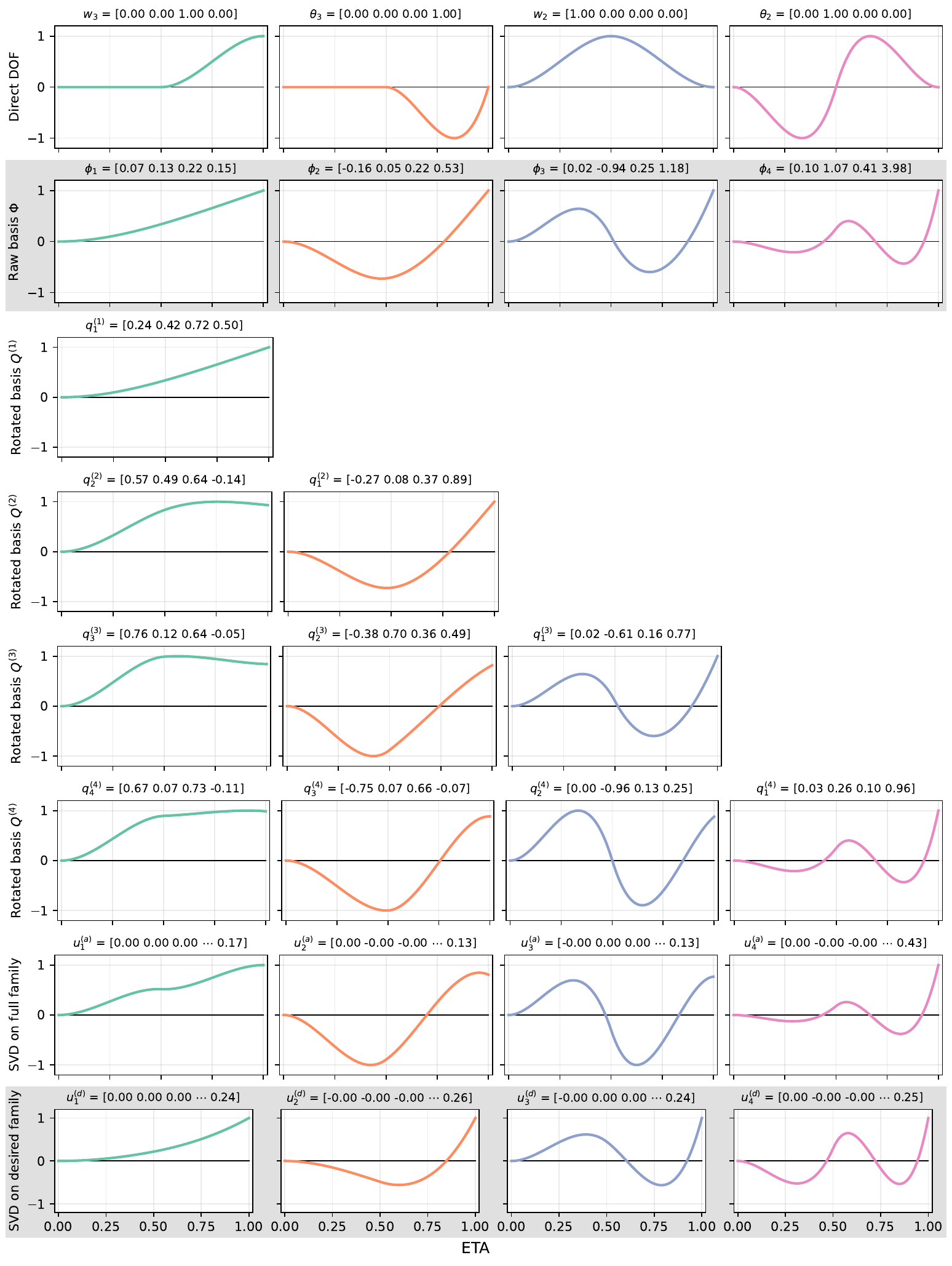}}

}

\caption{\label{fig-two-beam-comparison}Basis-direction comparison for
the 2-element beam. Rows show direct beam DOFs, raw MPM modes \(\Phi\),
selection-specific rotated bases \(Q^{(1)}\)-\(Q^{(4)}\), SVD trained on
the full arbitrary family, and SVD trained on the desired dihedral
family. Columns align comparable deflection shapes across methods.
Rotated vectors are reordered for this visual comparison, while their
\(q_i^{(s)}\) labels retain the QR-column indices.}

\end{figure}

Figure~\ref{fig-two-beam-comparison} compares the basis directions,
normalised by their peak absolute displacement. The direct beam-DOF row
shows isolated reduced-coordinate motions, while the raw and rotated MPM
rows show coupled directions derived from the pseudo-structure. The
rotated vectors are arranged by comparable shape across rows, while
their \(q_i^{(s)}\) labels retain the QR-column indices. The two SVD
rows then show what a data-driven reduction learns when trained on the
broader arbitrary family and on the desired dihedral family.

The selected bases were compared through a reconstruction study on the
full arbitrary and desired dihedral beam-shape libraries. For each
target family, the direct, raw MPM, and SVD directions were greedily
ranked by forward selection. At each selected count \(s\), the rotated
comparison factorised the corresponding raw MPM subset and used all
columns of \(Q^{(s)}\). The sampled shapes were reconstructed by least
squares using the selected basis at each count. A sample was counted as
covered when its relative reconstruction error remained below
\(\varepsilon = 0.01\). Table~\ref{tbl-beam-full-family} and
Table~\ref{tbl-beam-desired-family} report the resulting coverage counts
for one to four basis vectors.

\begin{longtable}[]{@{}llllll@{}}
\toprule\noalign{}
Method & Basis-vector order & 1 vector & 2 vectors & 3 vectors & 4
vectors \\
\midrule\noalign{}
\endfirsthead
\toprule\noalign{}
Method & Basis-vector order & 1 vector & 2 vectors & 3 vectors & 4
vectors \\
\midrule\noalign{}
\endhead
\bottomrule\noalign{}
\tabularnewline
\caption{Greedy forward-selection coverage counts for the full arbitrary
beam-shape family, reported for one to four basis vectors and counted
within
\(\varepsilon = 0.01\).}\label{tbl-beam-full-family}\tabularnewline
\endlastfoot
Direct DOF & \([w_3, w_2, \theta_2, \theta_3]\) & 0 (0\%) & 25 (0\%) &
925 (9\%) & 10000 (100\%) \\
Raw basis \(\Phi\) & \([\phi_1, \phi_2, \phi_3, \phi_4]\) & 1 (0\%) & 13
(0\%) & 786 (8\%) & 10000 (100\%) \\
Rotated basis \(Q^{(1)}\) & \([q^{(1)}_{1}]\) & 1 (0\%) & & & \\
Rotated basis \(Q^{(2)}\) & \([q^{(2)}_{1}, q^{(2)}_{2}]\) & & 13 (0\%)
& & \\
Rotated basis \(Q^{(3)}\) & \([q^{(3)}_{1}, q^{(3)}_{2}, q^{(3)}_{3}]\)
& & & 786 (8\%) & \\
Rotated basis \(Q^{(4)}\) &
\([q^{(4)}_{1}, q^{(4)}_{2}, q^{(4)}_{3}, q^{(4)}_{4}]\) & & & & 10000
(100\%) \\
SVD on full family & \([u^{(a)}_1, u^{(a)}_2, u^{(a)}_3, u^{(a)}_4]\) &
0 (0\%) & 25 (0\%) & 973 (10\%) & 10000 (100\%) \\
SVD on desired family & \([u^{(d)}_1, u^{(d)}_2, u^{(d)}_3, u^{(d)}_4]\)
& 0 (0\%) & 14 (0\%) & 416 (4\%) & 10000 (100\%) \\
\end{longtable}

\begin{longtable}[]{@{}llllll@{}}
\toprule\noalign{}
Method & Basis-vector order & 1 vector & 2 vectors & 3 vectors & 4
vectors \\
\midrule\noalign{}
\endfirsthead
\toprule\noalign{}
Method & Basis-vector order & 1 vector & 2 vectors & 3 vectors & 4
vectors \\
\midrule\noalign{}
\endhead
\bottomrule\noalign{}
\tabularnewline
\caption{Greedy forward-selection coverage counts for the desired
dihedral beam-shape family, reported for one to four basis vectors and
counted within
\(\varepsilon = 0.01\).}\label{tbl-beam-desired-family}\tabularnewline
\endlastfoot
Direct DOF & \([w_3, \theta_3, w_2, \theta_2]\) & 0 (0\%) & 14 (1\%) &
46 (5\%) & 1000 (100\%) \\
Raw basis \(\Phi\) & \([\phi_1, \phi_2, \phi_3, \phi_4]\) & 21 (2\%) &
196 (20\%) & 510 (51\%) & 1000 (100\%) \\
Rotated basis \(Q^{(1)}\) & \([q^{(1)}_{1}]\) & 21 (2\%) & & & \\
Rotated basis \(Q^{(2)}\) & \([q^{(2)}_{1}, q^{(2)}_{2}]\) & & 196
(20\%) & & \\
Rotated basis \(Q^{(3)}\) & \([q^{(3)}_{1}, q^{(3)}_{2}, q^{(3)}_{3}]\)
& & & 510 (51\%) & \\
Rotated basis \(Q^{(4)}\) &
\([q^{(4)}_{1}, q^{(4)}_{2}, q^{(4)}_{3}, q^{(4)}_{4}]\) & & & & 1000
(100\%) \\
SVD on full family & \([u^{(a)}_1, u^{(a)}_2, u^{(a)}_4, u^{(a)}_3]\) &
0 (0\%) & 0 (0\%) & 63 (6\%) & 1000 (100\%) \\
SVD on desired family & \([u^{(d)}_1, u^{(d)}_2, u^{(d)}_3, u^{(d)}_4]\)
& 4 (0\%) & 161 (16\%) & 638 (64\%) & 1000 (100\%) \\
\end{longtable}

The basis directions in Figure~\ref{fig-two-beam-comparison} distinguish
basis design from numerical conditioning. The raw MPM modes remain
smooth, coupled, and interpretable, whereas the rotated families
\(Q^{(s)}\) change visibly with selected count \(s\) because each
pivoted QR factorisation is built from a different selected subset
\(\Phi_s\). For each fixed selected count, \(Q^{(s)}\) spans the same
modal subspace as \(\Phi_s\), but pivoting changes the column order and
QR replaces the raw modes with Euclidean-orthonormal directions. The
resulting \(q_i^{(s)}\) vectors are therefore less stable to interpret
than the underlying pseudo-structural modes. The SVD directions were
also strongly library-dependent. The full arbitrary family included
sagging and non-monotone behaviour, while the desired dihedral family
shifted the dominant vectors toward the same bending patterns seen in
the MPM modes.

On the desired dihedral family, the raw two-vector MPM basis covered 196
of 1,000 sampled shapes within the reconstruction tolerance, compared
with 161 of 1,000 for the SVD basis trained directly on that same
family. MPM performed worse on the broader arbitrary family, but this
was acceptable because the full family included shape features such as
sagging and overshoot that would be unlikely in transport-aircraft
dihedral curves. The principal-angle comparison followed the same trend:
the principal-angle vectors between the raw MPM and desired-family SVD
selected spaces were \([7.5^\circ]\), \([0.2^\circ, 21.4^\circ]\), and
\([0.0^\circ, 0.0^\circ, 21.8^\circ]\) for one, two, and three
directions, respectively, whereas the third principal angle against the
full-family SVD basis increased to \(72.4^\circ\). This reflected a
benefit of MPM, where a tuned pseudo-structure can generate useful
admissible shapes similar to SVD but without requiring a prior library
of desired shapes.

\subsection{Basis-Design Tuning: 3D Thin Airfoil
Wing}\label{sec-3d-wing-basis-design}

\begin{figure}

\centering{

\includegraphics[width=0.75\linewidth,height=\textheight,keepaspectratio]{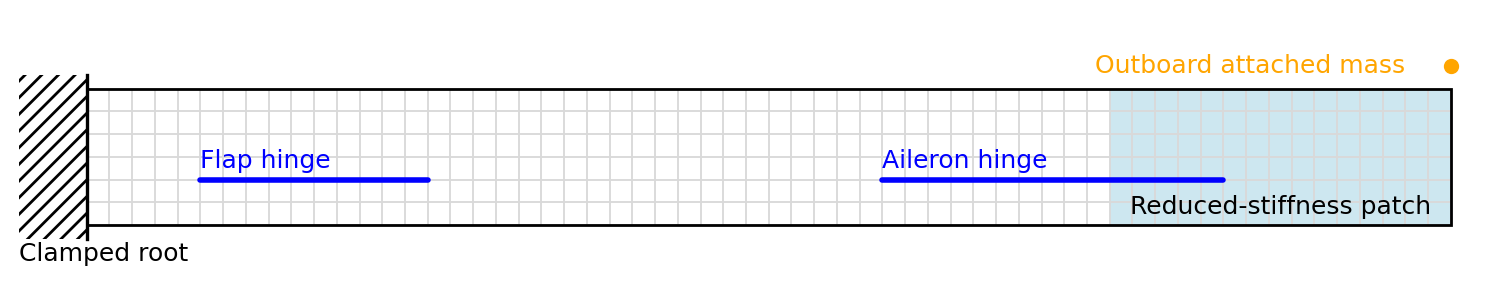}

}

\caption{\label{fig-thin-wing-setup}Rectangular cantilever thin-wing
used for the basis-design study. The same geometry underlies the
baseline, stiffness-tuned, mass-tuned, protected-hinge, and composite
cases.}

\end{figure}

A simple rectangular cantilevered shell wing, shown in
Figure~\ref{fig-thin-wing-setup}, was used to demonstrate how MPM basis
functions can be tuned on a 3D planform. The structural modes for the
cases shown in Table~\ref{tbl-thin-wing-setup} were generated with
NASTRAN SOL 103, and the deformed geometry families were evaluated using
classical lifting line theory {[}32{]}. This study builds on the 2D
two-element beam model in Section~\ref{sec-cantilever-beam} by showing
how basis design through stiffness edits, added point masses, and
constraints influences the admissible shapes. Each tuned case was
compared against a baseline with uniform stiffness, a clamped root, no
added mass, and no hinge constraints at matched \(C_L = 0.4\) for
\(0^\circ \leq \alpha \leq 10^\circ\).

\begin{longtable}[]{@{}lccccccccc@{}}
\toprule\noalign{}
Case & \(c\) {[}m{]} & \(b/2\) {[}m{]} & \(t\) {[}m{]} & \(E_0\)
{[}GPa{]} & \(\nu\) & \(\rho\) {[}kg m\(^{-3}\){]} & \(E_p/E_0\) &
\(M_p\) {[}kg{]} & Hinge \\
\midrule\noalign{}
\endfirsthead
\toprule\noalign{}
Case & \(c\) {[}m{]} & \(b/2\) {[}m{]} & \(t\) {[}m{]} & \(E_0\)
{[}GPa{]} & \(\nu\) & \(\rho\) {[}kg m\(^{-3}\){]} & \(E_p/E_0\) &
\(M_p\) {[}kg{]} & Hinge \\
\midrule\noalign{}
\endhead
\bottomrule\noalign{}
\tabularnewline
\caption[]{Thin-wing basis-design cases. Geometry and base material
remain fixed across all cases, while the final columns identify the
stiffness, added-mass, and hinge-line edits introduced in each case. The
Hinge column reports the number of multi-point constraints
(MPCs).}\label{tbl-thin-wing-setup}\tabularnewline
\endlastfoot
Baseline (untuned) & \multirow{5}{*}{\centering\arraybackslash  6.0} &
\multirow{5}{*}{\centering\arraybackslash  60.0} &
\multirow{5}{*}{\centering\arraybackslash  0.1} &
\multirow{5}{*}{\centering\arraybackslash  100} &
\multirow{5}{*}{\centering\arraybackslash  0.30} &
\multirow{5}{*}{\centering\arraybackslash  2000} & N/A & N/A & N/A \\
Stiffness tuning & & & & & & & 0.20 & N/A & N/A \\
Mass tuning & & & & & & & N/A & 900 & N/A \\
Protected hinge lines & & & & & & & N/A & N/A & 2 MPCs \\
Composite case & & & & & & & 0.20 & 900 & 2 MPCs \\
\end{longtable}

\begin{figure}

\centering{

\includegraphics[width=0.88\linewidth,height=\textheight,keepaspectratio]{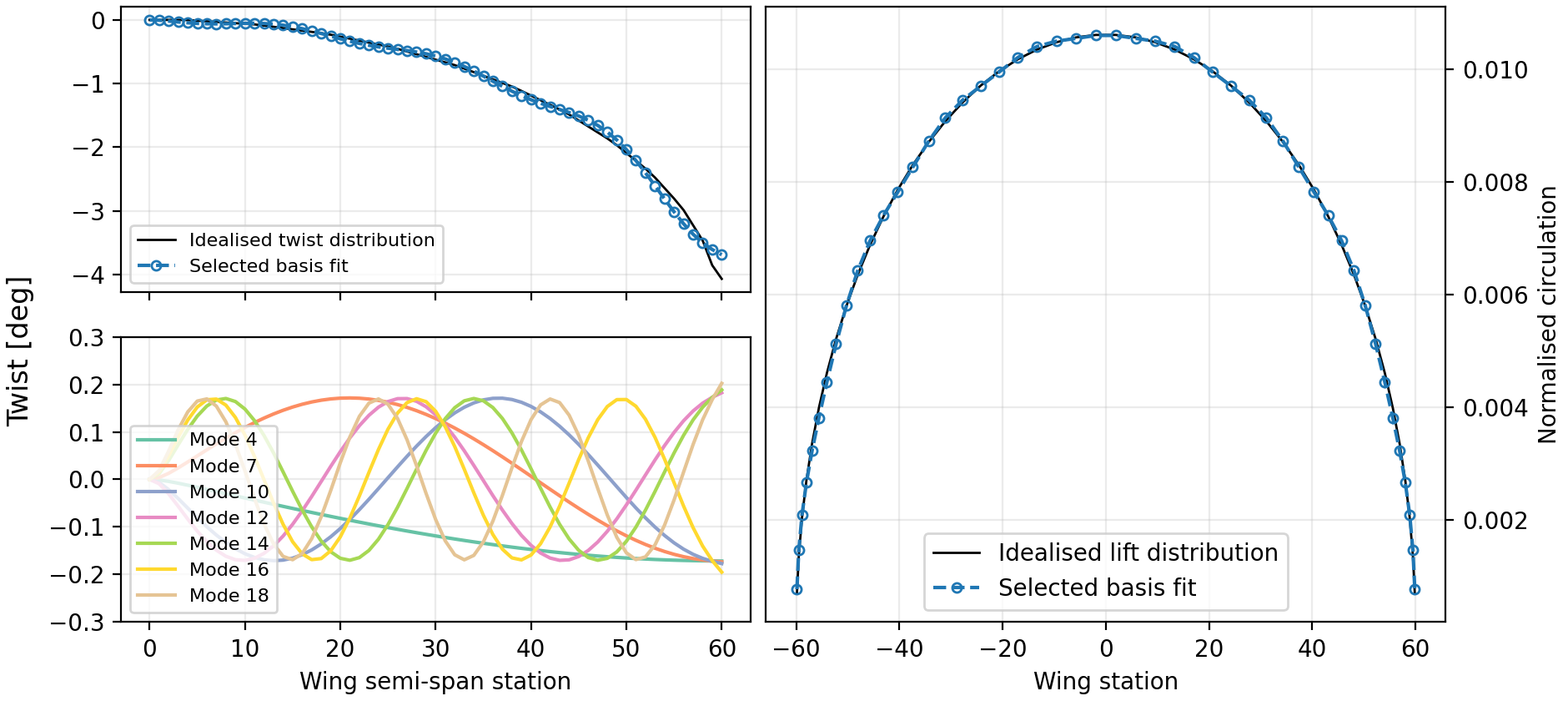}

}

\caption{\label{fig-thin-wing-baseline-fit-summary}Untuned baseline
thin-wing fit to the elliptic lift-distribution target at \(C_L = 0.4\).
Seven modes are selected, and the resulting twist remains visibly
noisy.}

\end{figure}

The untuned baseline was first evaluated using the first 50 structural
modes and a greedy forward-selection scheme to choose basis functions
that would maximise the Prandtl span efficiency factor. While the target
spanwise efficiency factor \(e \geq 0.999\) was achieved at
\(C_L = 0.4\), Figure~\ref{fig-thin-wing-baseline-fit-summary} shows
that seven modes (4, 7, 10, 12, 14, 16, 18) were required and the wing
twist was visibly noisy throughout the span compared to the ideal twist
distribution for an elliptic lift distribution.

\begin{figure}

\centering{

\includegraphics[width=0.88\linewidth,height=\textheight,keepaspectratio]{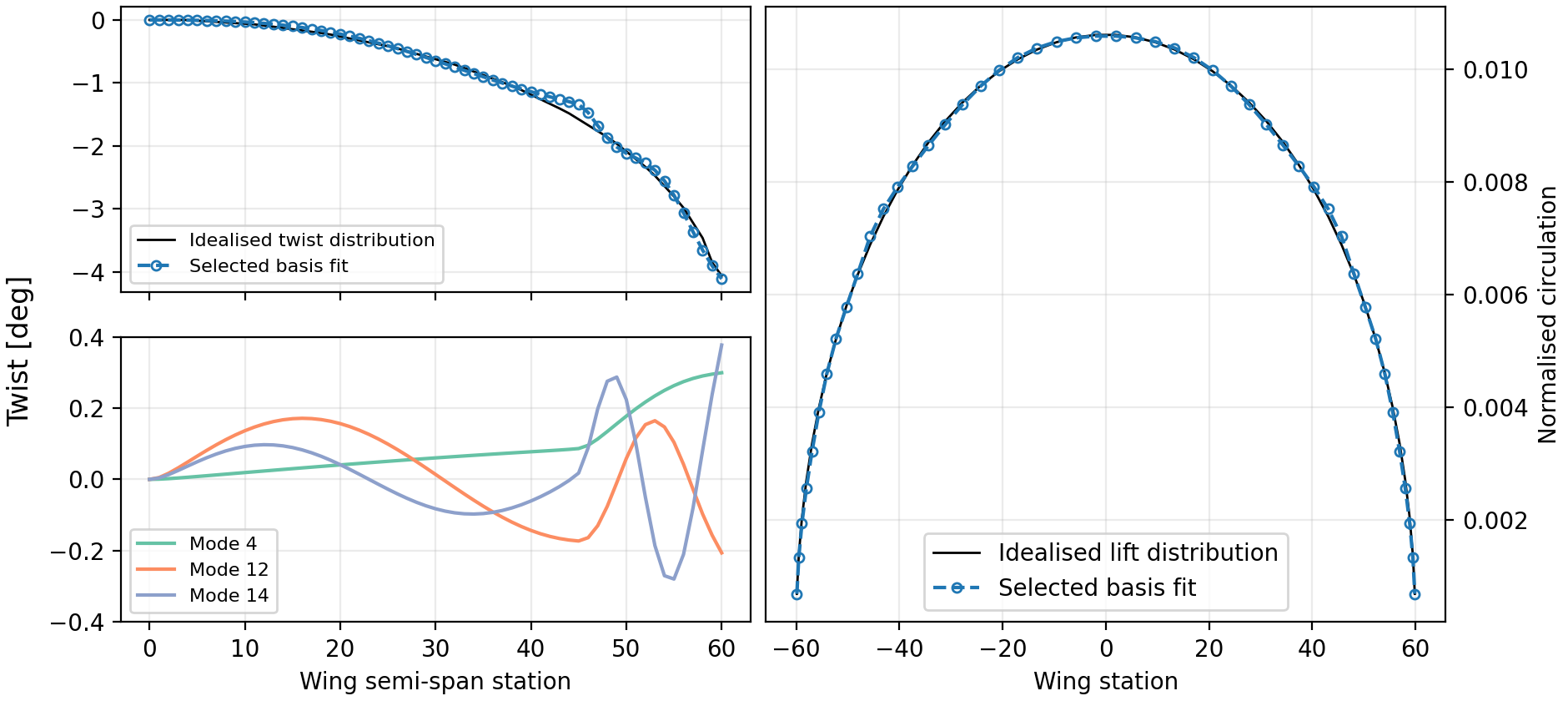}

}

\caption{\label{fig-thin-wing-stiffness-fit-summary}Stiffness-tuned
thin-wing fit for the elliptic target. The outboard soft patch reduces
the selected set to three modes but introduces a visible kink near the
stiffness transition.}

\end{figure}

In comparison, the stiffness-tuning case with a full-chord outboard soft
patch reduced the number of selected modes required, with only three
modes required (4, 12, 14) to reach \(e \geq 0.999\). The twist was
smoother in the inboard and outboard regions than in the baseline case,
but the stiffness transition at \(y = 45\) introduced a visible kink, as
shown in Figure~\ref{fig-thin-wing-stiffness-fit-summary}. For more
complex geometries, patch number, size, location, and stiffness all
become tuning variables. This can improve family alignment, but it also
increases modelling effort and implementation burden.

\begin{figure}

\centering{

\includegraphics[width=0.88\linewidth,height=\textheight,keepaspectratio]{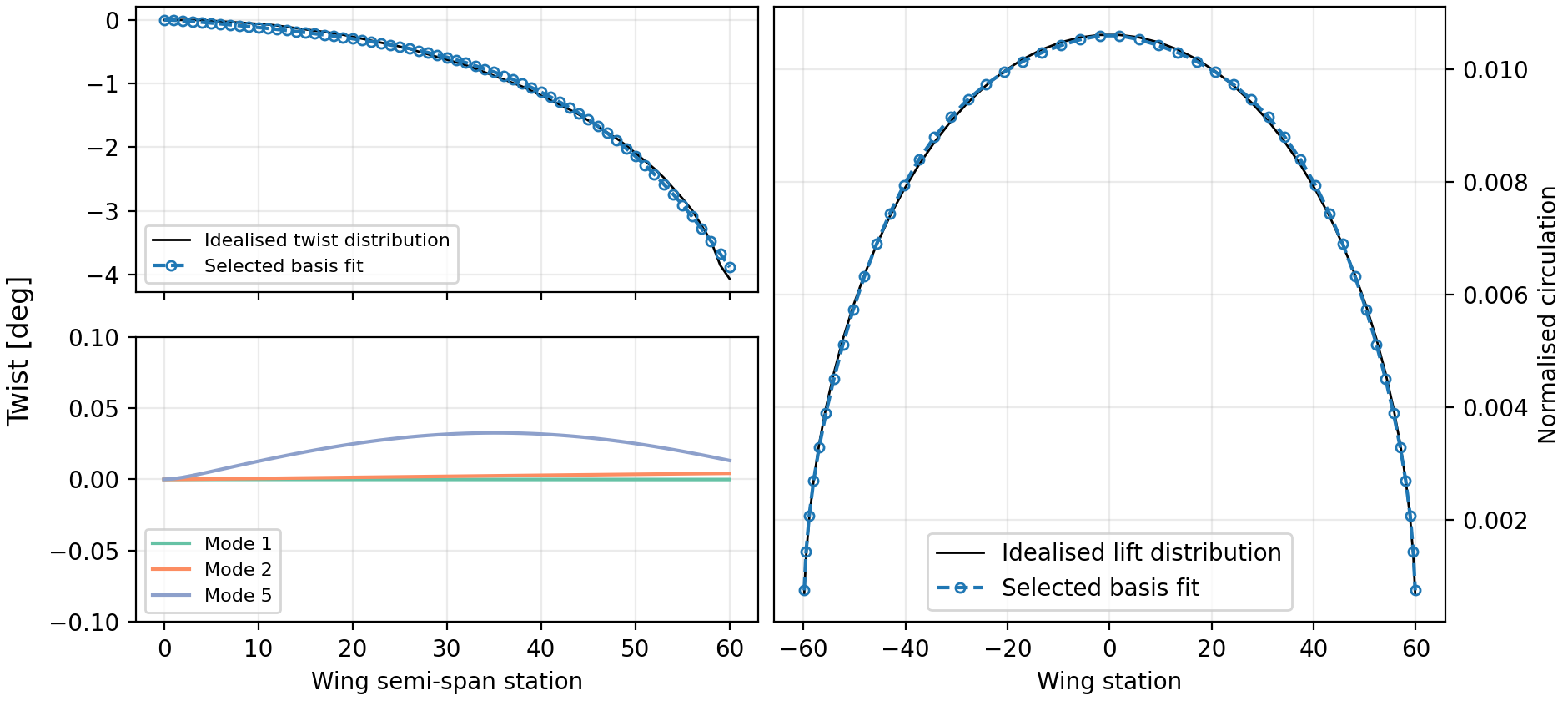}

}

\caption{\label{fig-thin-wing-mass-fit-summary}Mass-tuned thin-wing fit
for the elliptic target. The attached tip mass also reaches the target
with three modes, while giving a smoother twist distribution than the
stiffness-tuned case.}

\end{figure}

In practice, mass tuning is a more effective basis-design tool. For the
mass-tuning case, a single NASTRAN CONM2 point mass was connected to the
wing tip with NASTRAN RBE3s. Similar to stiffness tuning, the
mass-tuning case also reduced the number of modes required to reach
\(e \geq 0.999\) to only three modes (1, 2, 5), but it produced a much
smoother twist distribution as shown in
Figure~\ref{fig-thin-wing-mass-fit-summary}. Only the general location
of the point mass needed to be defined, in this example at the wing tip,
and varying the mass magnitude became the main tuning variable, which
was simpler to adjust than a distributed stiffness patch.

Increasing the attached mass shifts the selected modal family from a
global to a more localised outboard response, with increasing localised
washout at the wing tip. This is another benefit of MPM, where global
and local refinement can be achieved effectively through structural
parameter changes. In comparison, many methods in
Section~\ref{sec-related-work} require the reparameterisation of the
geometry.

\begin{figure}

\centering{

\includegraphics[width=0.88\linewidth,height=\textheight,keepaspectratio]{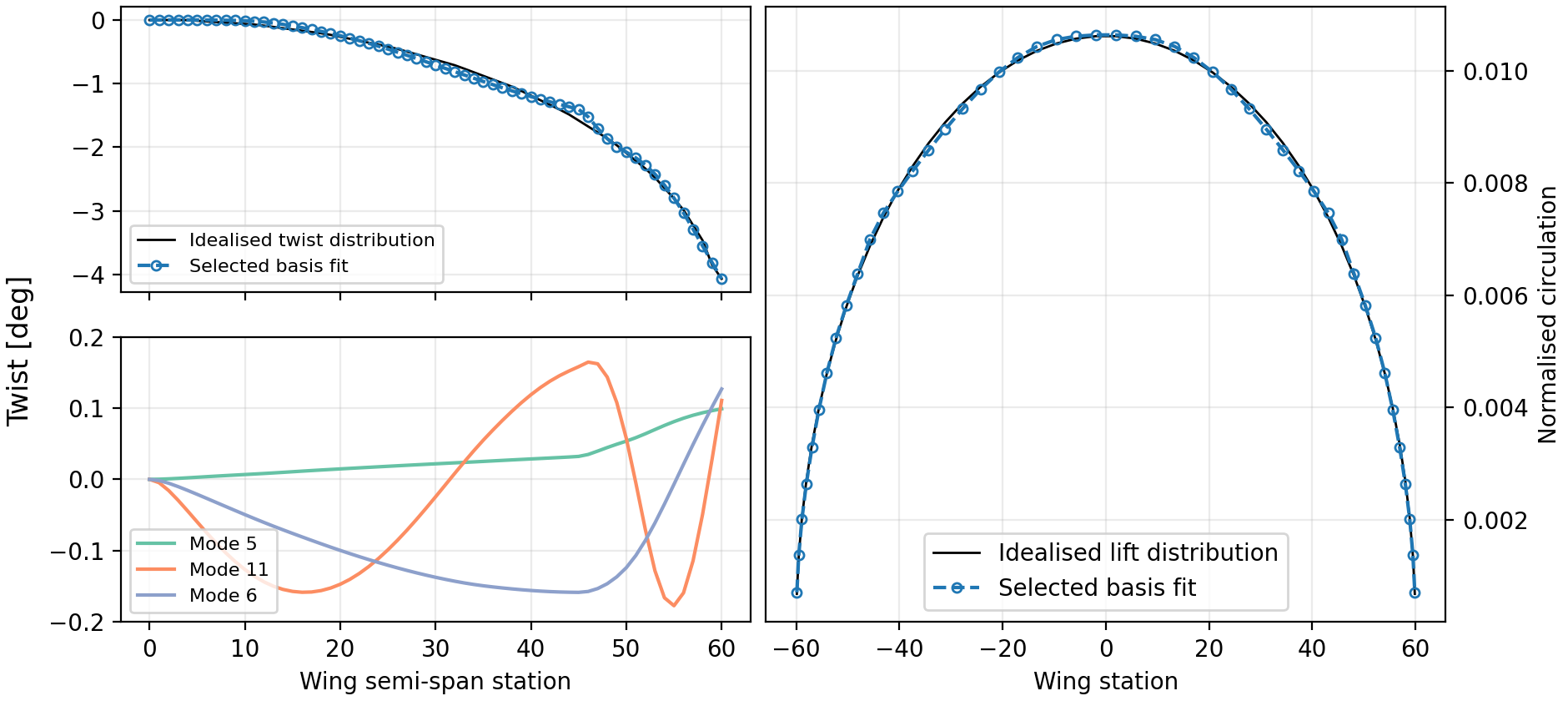}

}

\caption{\label{fig-thin-wing-composite-fit-summary}Composite thin-wing
case combining stiffness tuning, attached mass, and protected hinge
lines. The elliptic target is still met with three selected modes
despite the added admissibility constraints.}

\end{figure}

The composite pseudo-structure combined the outboard stiffness
reduction, attached mass, and protected hinge lines. Three selected
modes, containing both low- and higher-order content, achieved a span
efficiency of \(e \geq 0.999\).
Figure~\ref{fig-thin-wing-composite-fit-summary} shows the resulting
elliptic-target fit, indicating that basis-shaping edits and hinge-line
constraints can be combined while retaining a compact selected basis.

\begin{figure}

\begin{minipage}[t]{0.33\linewidth}

\centering{

\includegraphics[width=1\linewidth,height=\textheight,keepaspectratio]{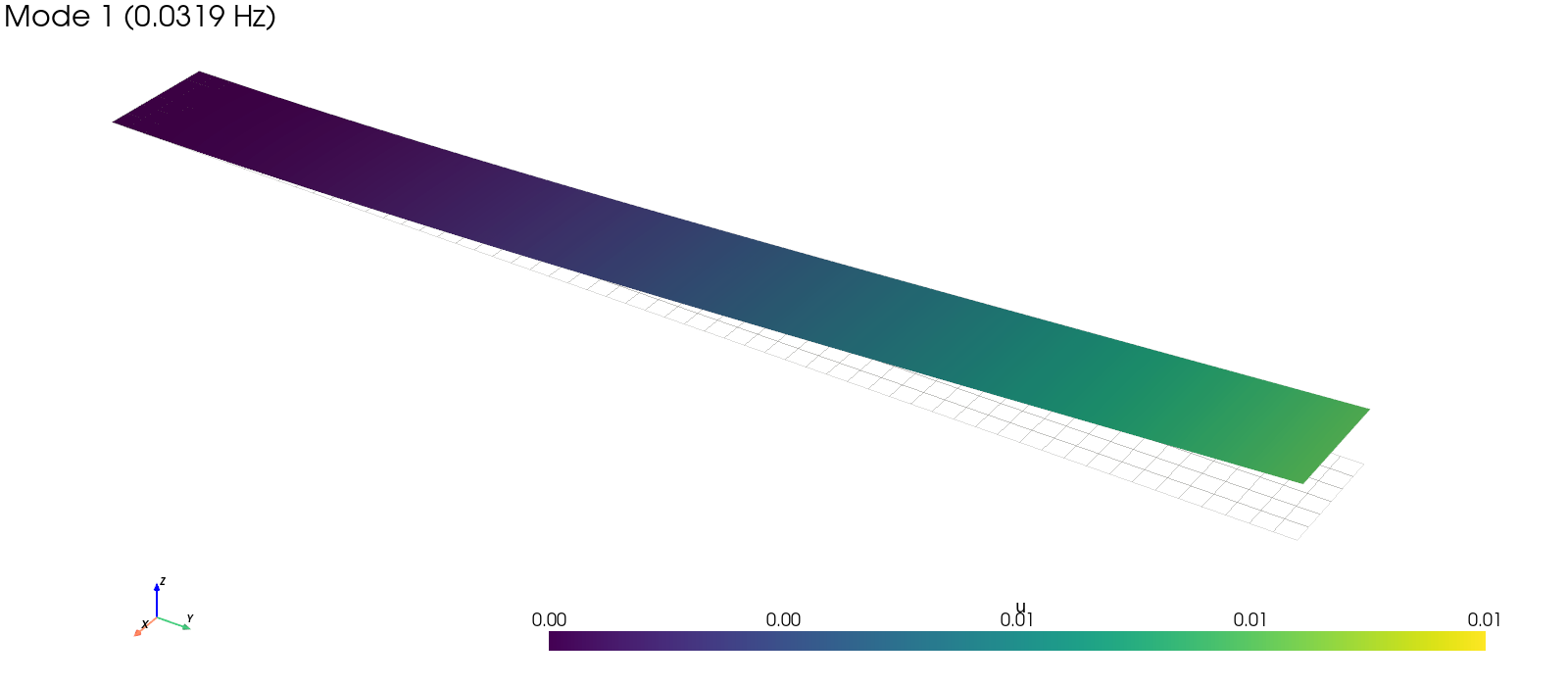}

}

\subcaption{\label{fig-base-mode-1}Baseline case, mode 1.}

\end{minipage}%
\begin{minipage}[t]{0.33\linewidth}

\centering{

\includegraphics[width=1\linewidth,height=\textheight,keepaspectratio]{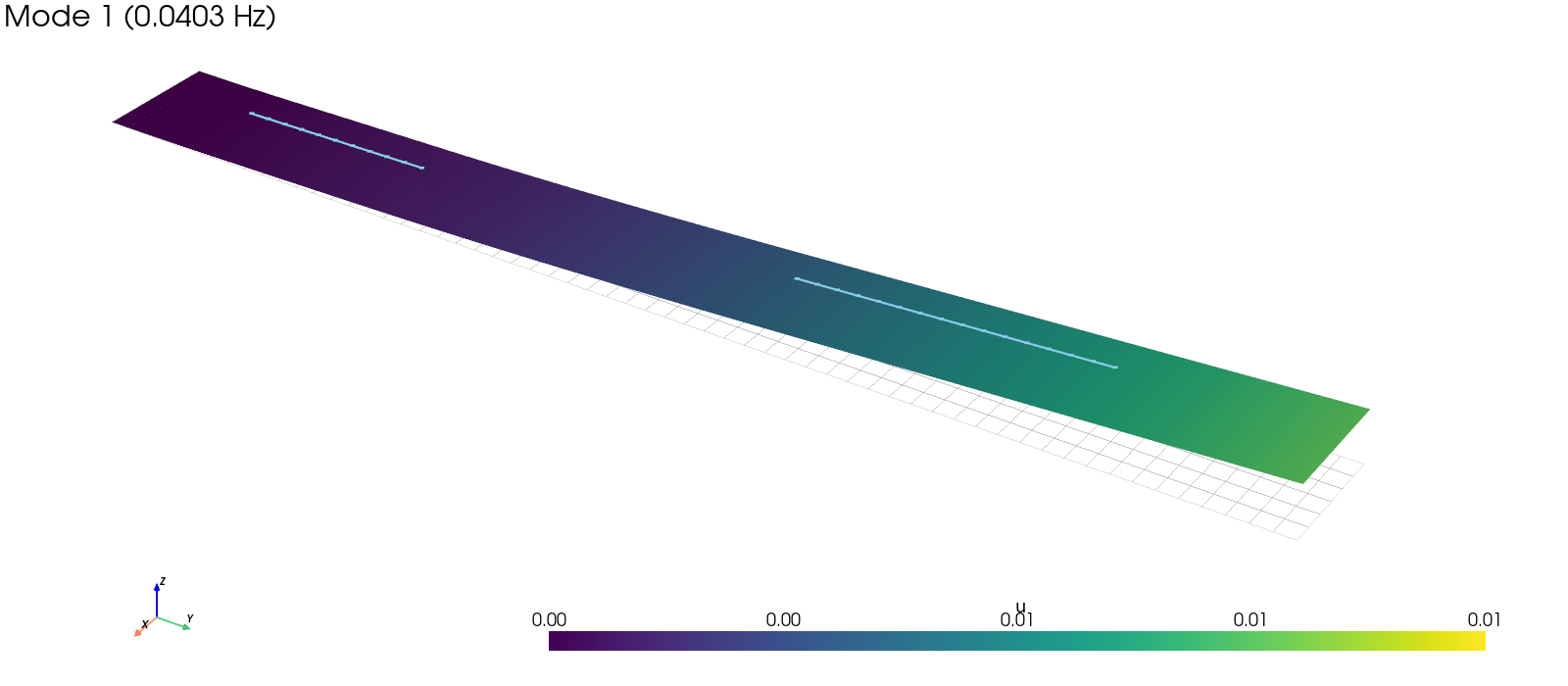}

}

\subcaption{\label{fig-hinge-mode-1}Protected-hinge case, mode 1.}

\end{minipage}%
\begin{minipage}[t]{0.33\linewidth}

\centering{

\includegraphics[width=1\linewidth,height=\textheight,keepaspectratio]{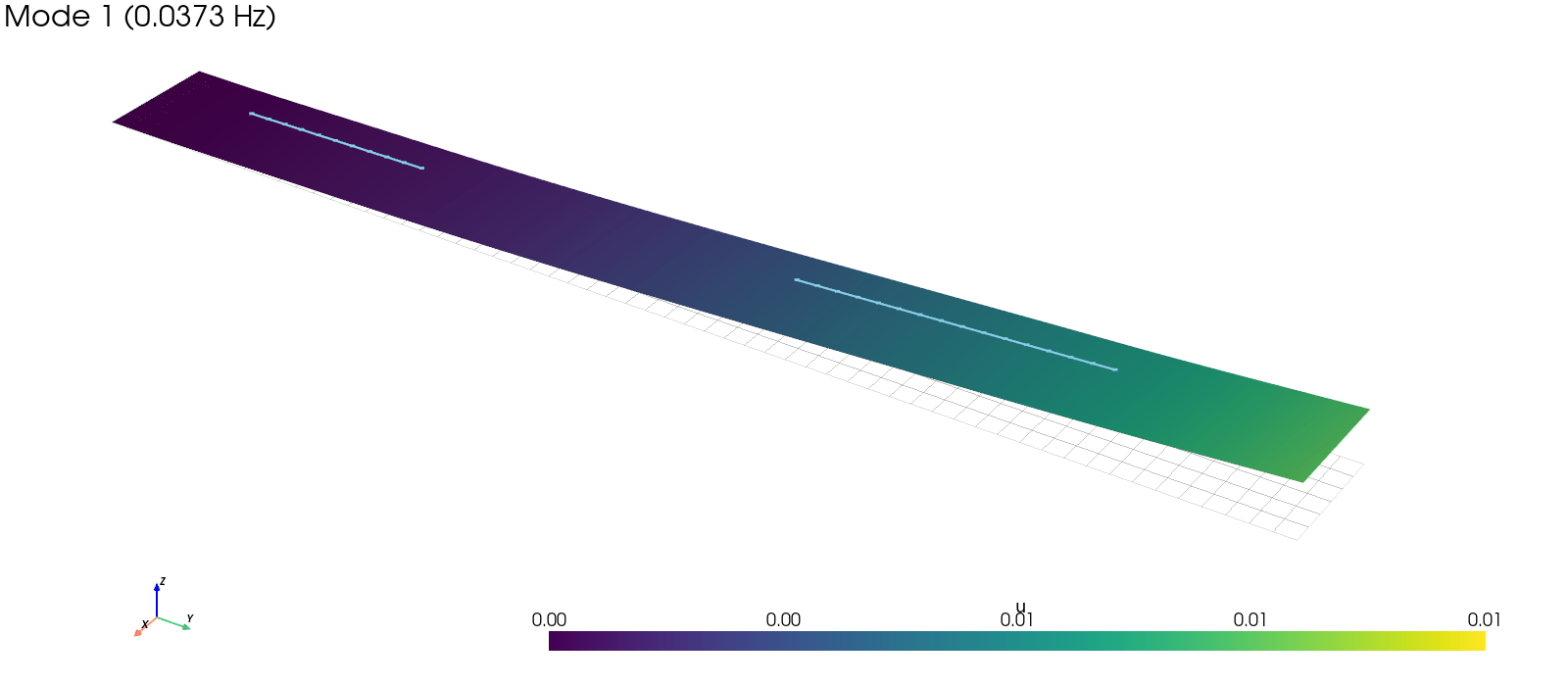}

}

\subcaption{\label{fig-composite-mode-1}Composite case, mode 1.}

\end{minipage}%
\newline
\begin{minipage}[t]{0.33\linewidth}

\centering{

\includegraphics[width=1\linewidth,height=\textheight,keepaspectratio]{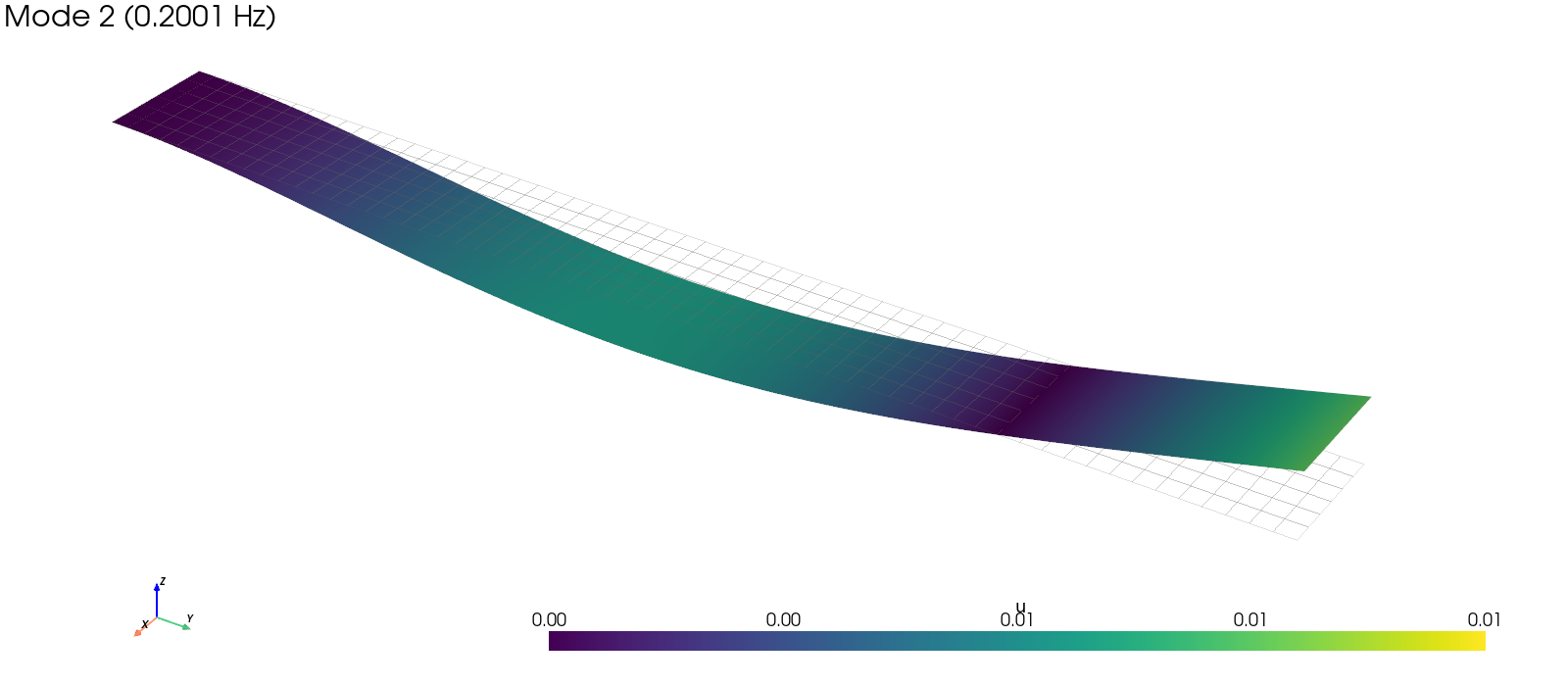}

}

\subcaption{\label{fig-base-mode-2}Baseline case, mode 2.}

\end{minipage}%
\begin{minipage}[t]{0.33\linewidth}

\centering{

\includegraphics[width=1\linewidth,height=\textheight,keepaspectratio]{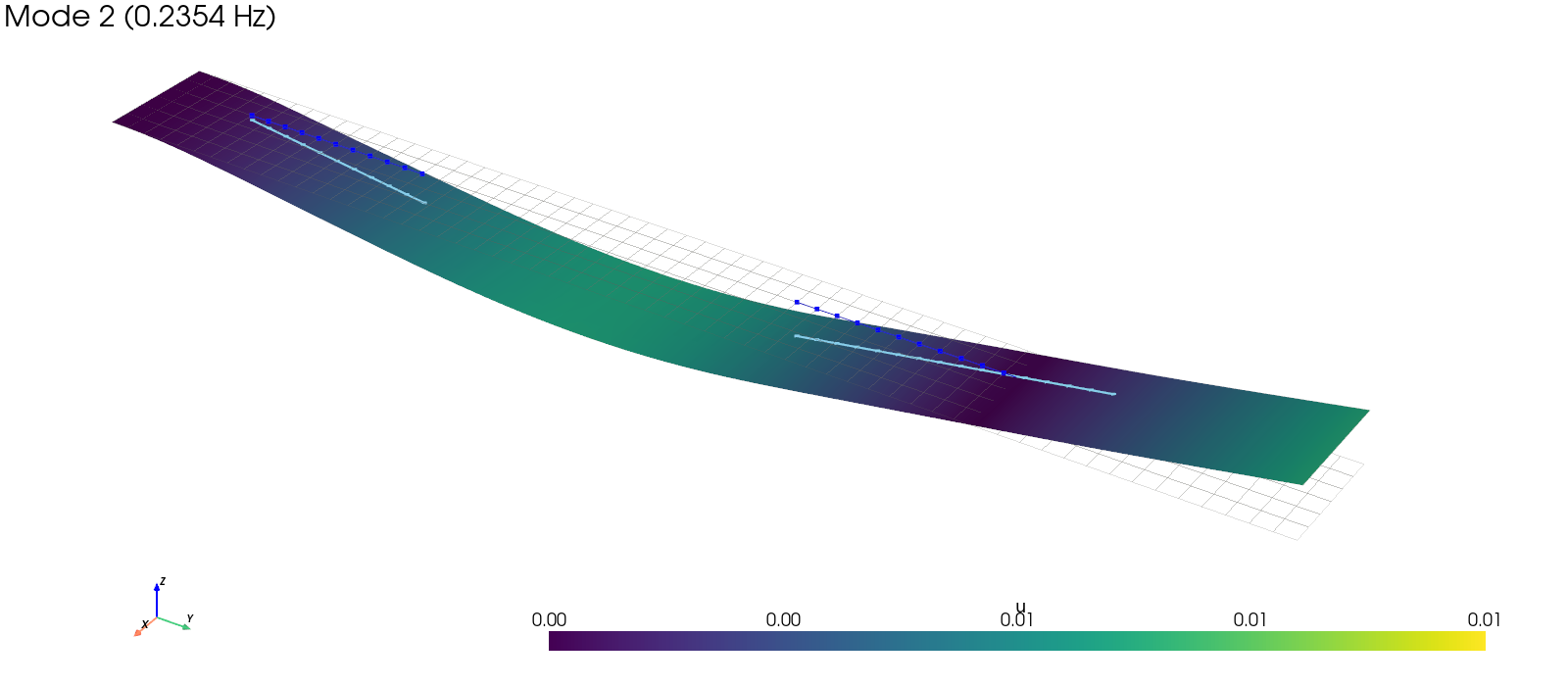}

}

\subcaption{\label{fig-hinge-mode-2}Protected-hinge case, mode 2.}

\end{minipage}%
\begin{minipage}[t]{0.33\linewidth}

\centering{

\includegraphics[width=1\linewidth,height=\textheight,keepaspectratio]{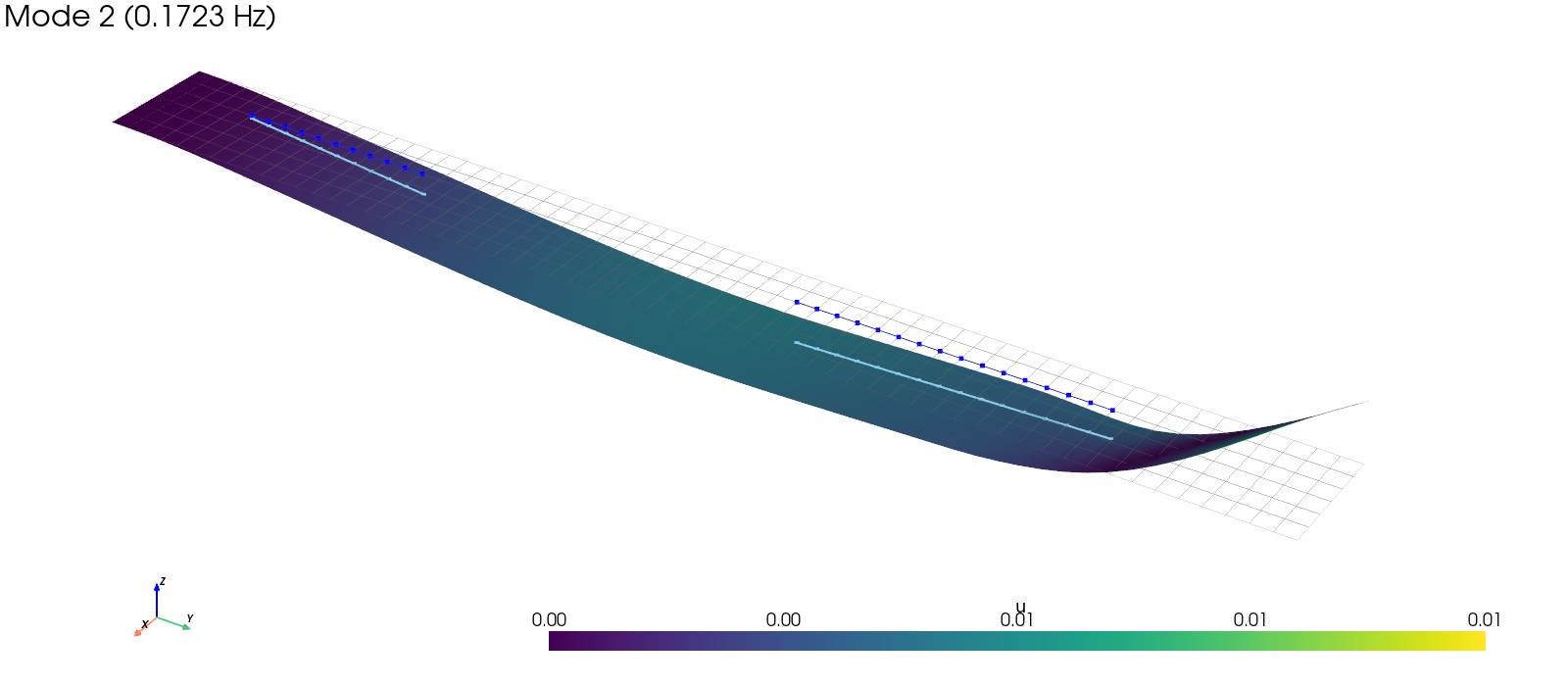}

}

\subcaption{\label{fig-composite-mode-2}Composite case, mode 2.}

\end{minipage}%
\newline
\begin{minipage}[t]{0.33\linewidth}

\centering{

\includegraphics[width=1\linewidth,height=\textheight,keepaspectratio]{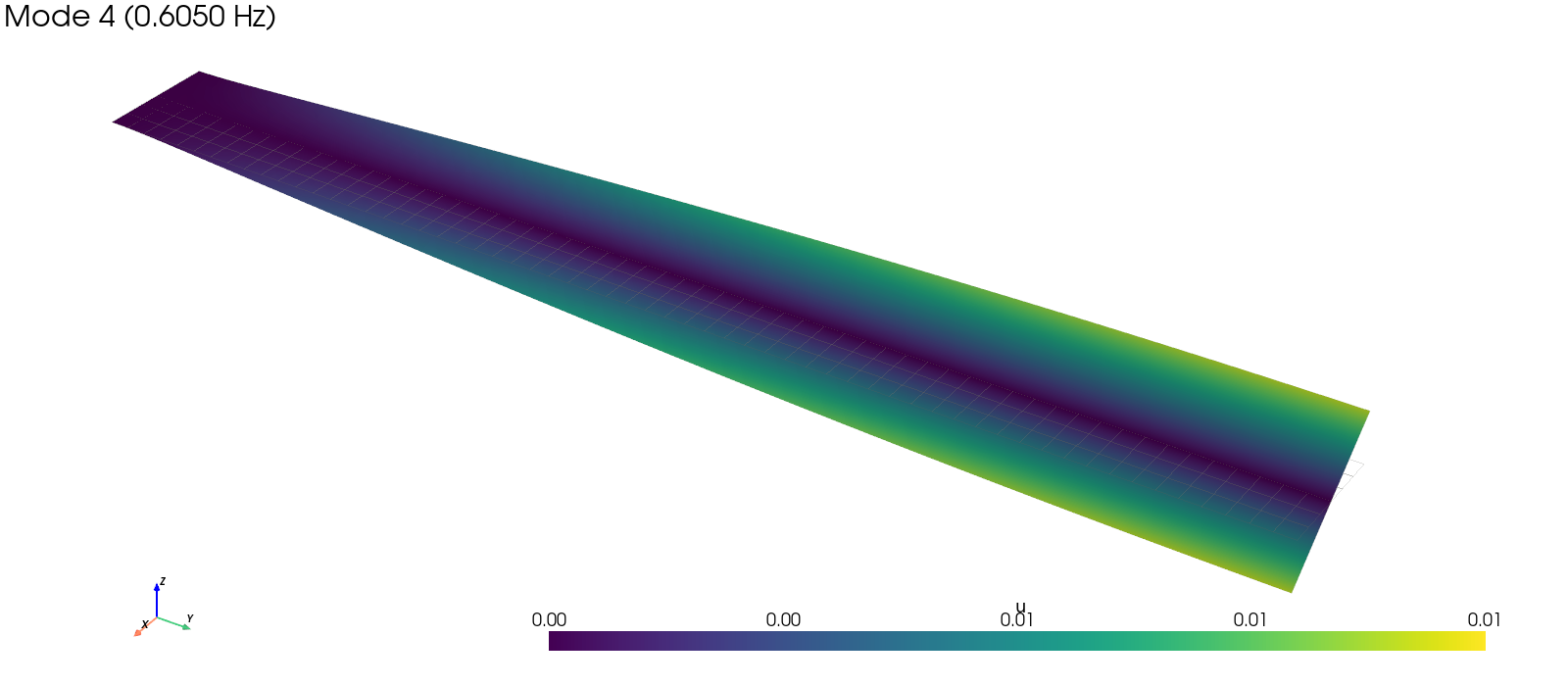}

}

\subcaption{\label{fig-base-mode-4}Baseline case, mode 4.}

\end{minipage}%
\begin{minipage}[t]{0.33\linewidth}

\centering{

\includegraphics[width=1\linewidth,height=\textheight,keepaspectratio]{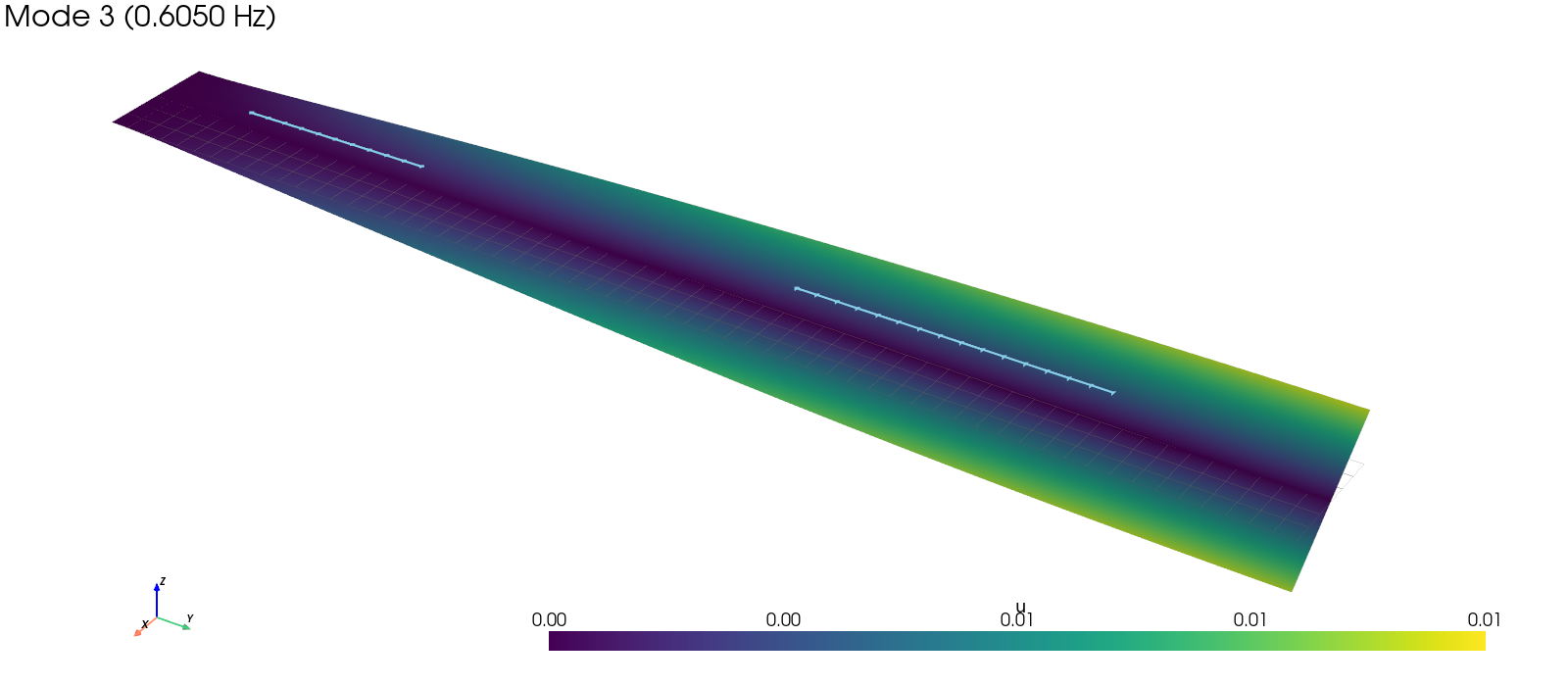}

}

\subcaption{\label{fig-hinge-mode-3}Protected-hinge case, mode 3.}

\end{minipage}%
\begin{minipage}[t]{0.33\linewidth}

\centering{

\includegraphics[width=1\linewidth,height=\textheight,keepaspectratio]{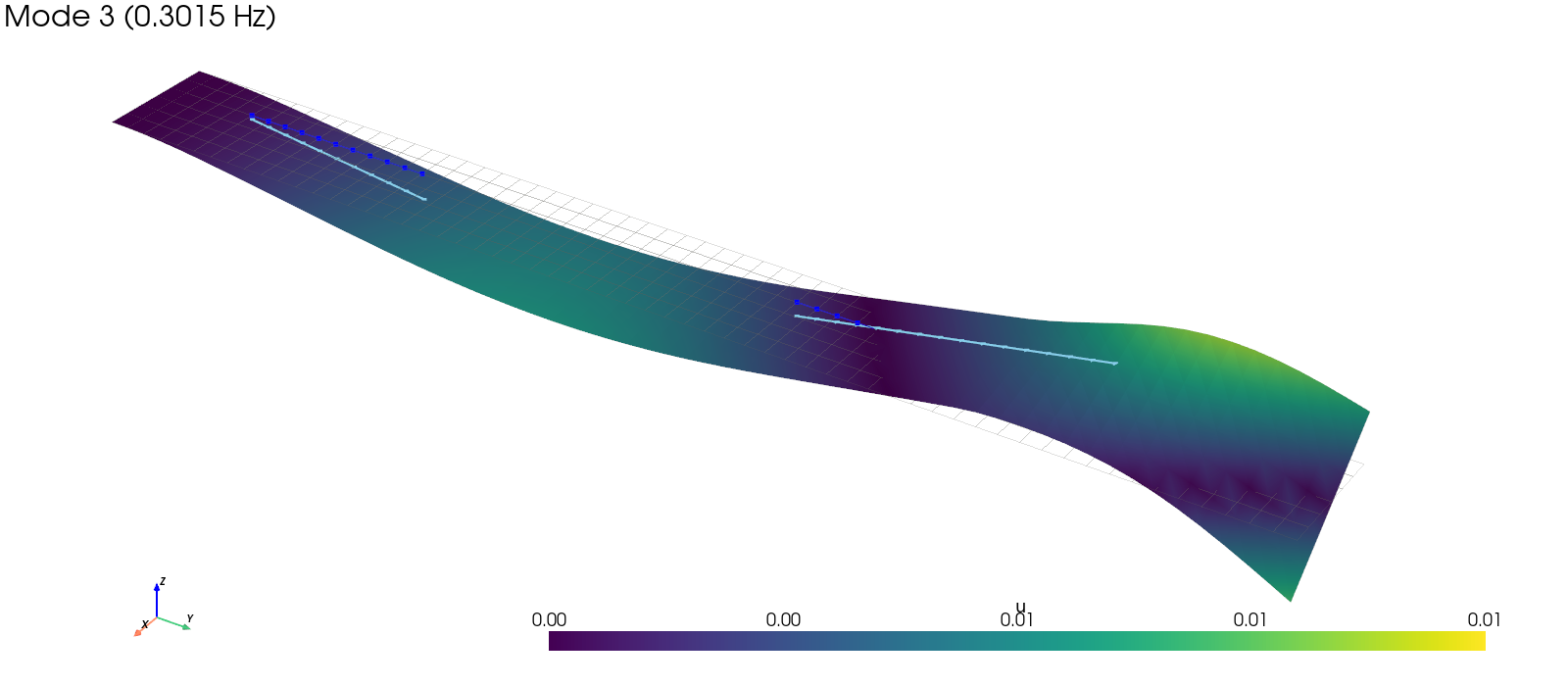}

}

\subcaption{\label{fig-composite-mode-3}Composite case, mode 3.}

\end{minipage}%
\newline
\begin{minipage}[t]{0.33\linewidth}

\centering{

\includegraphics[width=1\linewidth,height=\textheight,keepaspectratio]{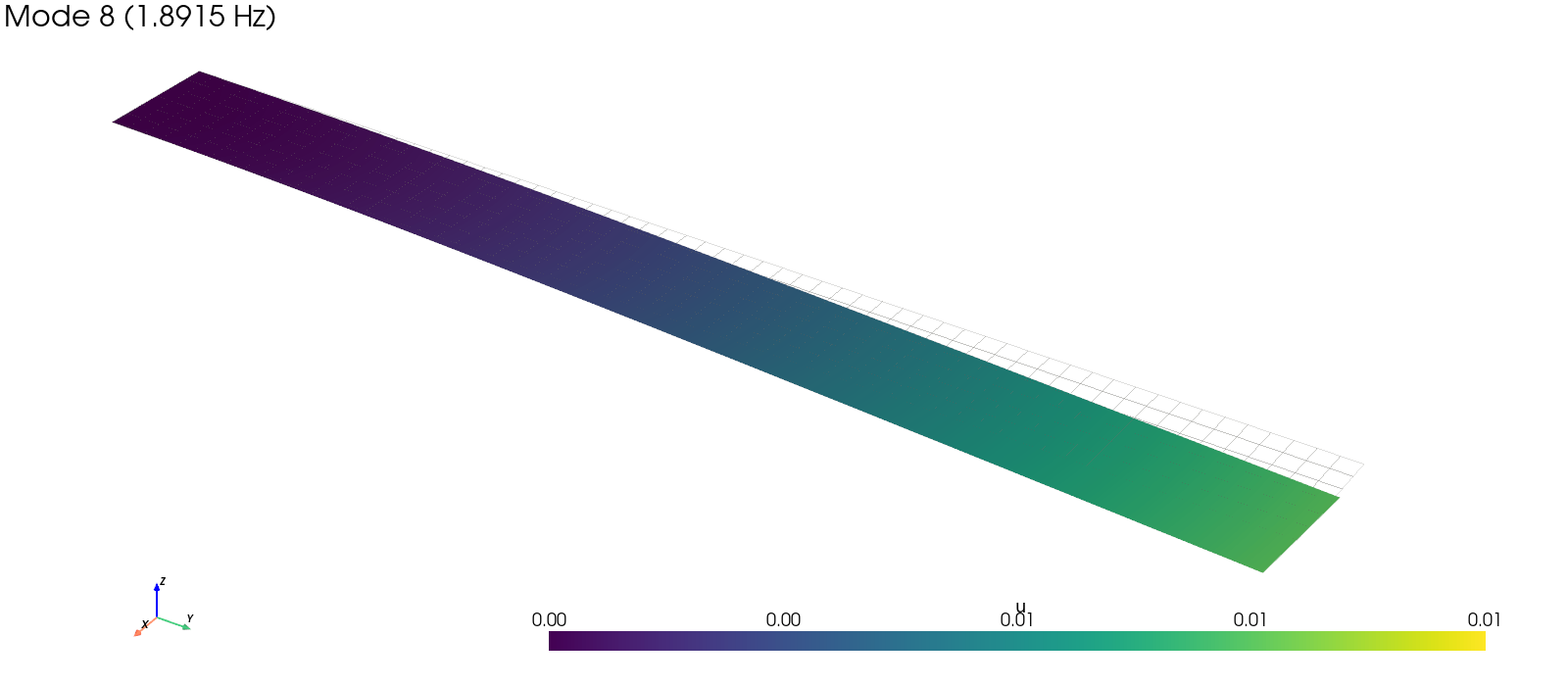}

}

\subcaption{\label{fig-base-mode-8}Baseline case, mode 8.}

\end{minipage}%
\begin{minipage}[t]{0.33\linewidth}

\centering{

\includegraphics[width=1\linewidth,height=\textheight,keepaspectratio]{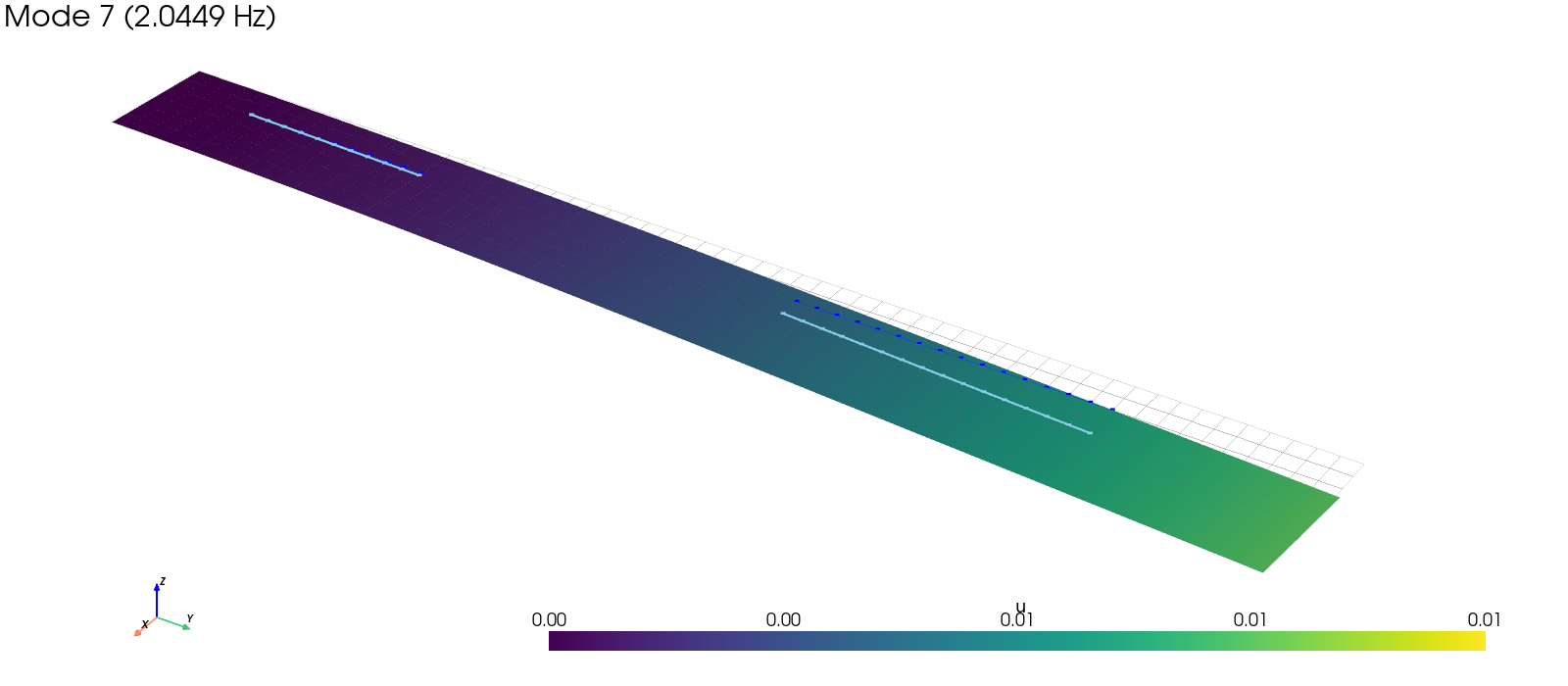}

}

\subcaption{\label{fig-hinge-mode-7}Protected-hinge case, mode 7.}

\end{minipage}%
\begin{minipage}[t]{0.33\linewidth}

\centering{

\includegraphics[width=1\linewidth,height=\textheight,keepaspectratio]{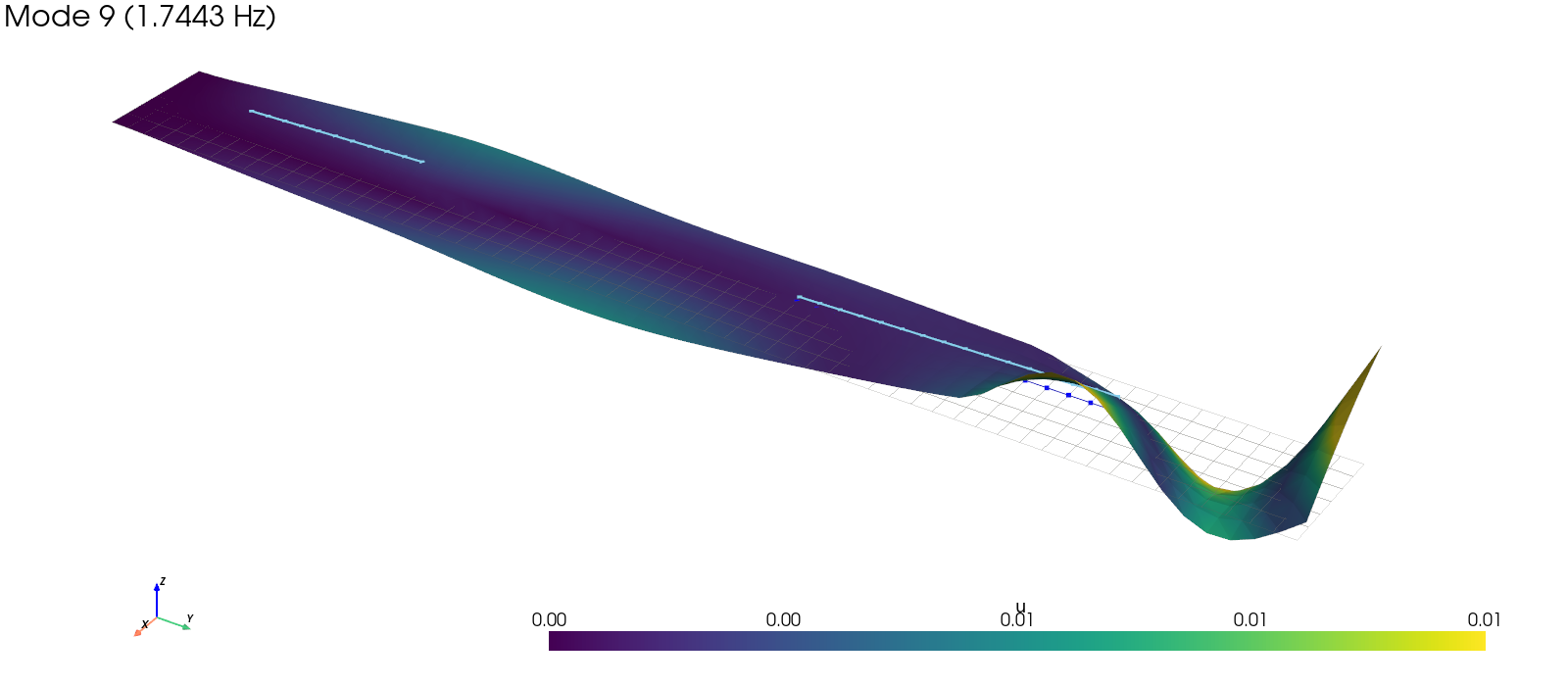}

}

\subcaption{\label{fig-composite-mode-9}Composite case, mode 9.}

\end{minipage}%

\caption{\label{fig-thin-wing-mode-comparison}Mode-shape comparison for
the baseline, protected-hinge, and composite thin-wing
pseudo-structures. Columns show baseline, protected-hinge, and composite
cases. Rows show modal-assurance-criterion (MAC)-matched first and
second out-of-plane bending, first torsion, and first in-plane bending
modes.}

\end{figure}

This also demonstrates another advantage of MPM, that it can represent
awkward geometric constraints directly in the admissible family.
Straight hinge lines may be required for manufacturing simplicity or for
morphing-wing concepts using continuous rotary actuators. In the
pseudo-structure, two NASTRAN multi-point constraints (MPCs) allowed
rigid-body translation and rotation of the hinge lines while still
allowing leading-edge and trailing-edge deformation. Alternative
parameterisations would typically need a bespoke analytical basis, a
customised FFD implementation, or a dedicated geometry library for
reduction. Selected mode shapes for the hinge lines are shown in
Figure~\ref{fig-thin-wing-mode-comparison}, including out-of-plane
bending modes associated with dihedral, torsional modes associated with
twist, and in-plane bending modes associated with sweep.

\subsection{Constraint Editing: Military Engine Intake Research Duct
(MEIRD)}\label{sec-meird}

\begin{figure}

\begin{minipage}[t]{0.50\linewidth}

\centering{

\includegraphics[width=0.55\linewidth,height=\textheight,keepaspectratio]{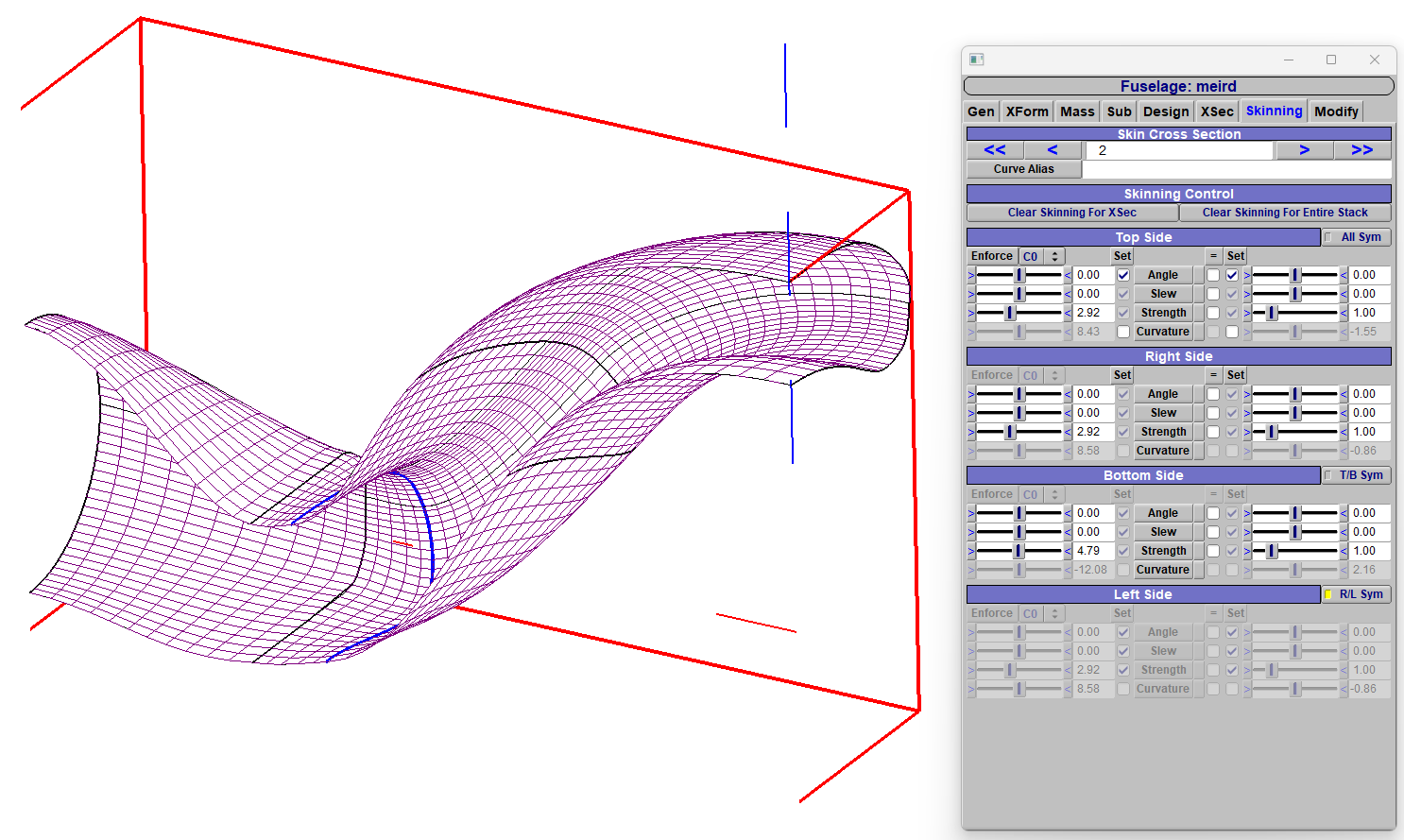}

}

\subcaption{\label{fig-meird-openvsp-param}Compact analytic: OpenVSP
parameterisation.}

\end{minipage}%
\begin{minipage}[t]{0.50\linewidth}

\centering{

\includegraphics[width=0.55\linewidth,height=\textheight,keepaspectratio]{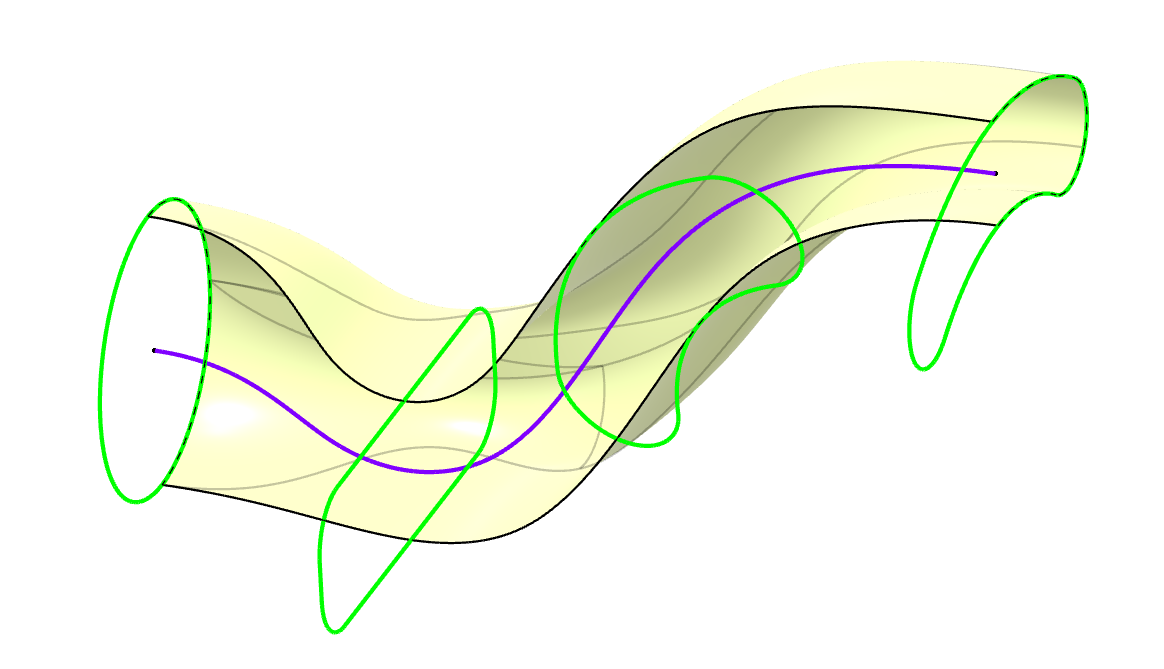}

}

\subcaption{\label{fig-meird-3dx-param}Direct geometry parameterisation:
CATIA B-rep.}

\end{minipage}%
\newline
\begin{minipage}[t]{0.50\linewidth}

\centering{

\includegraphics[width=0.55\linewidth,height=\textheight,keepaspectratio]{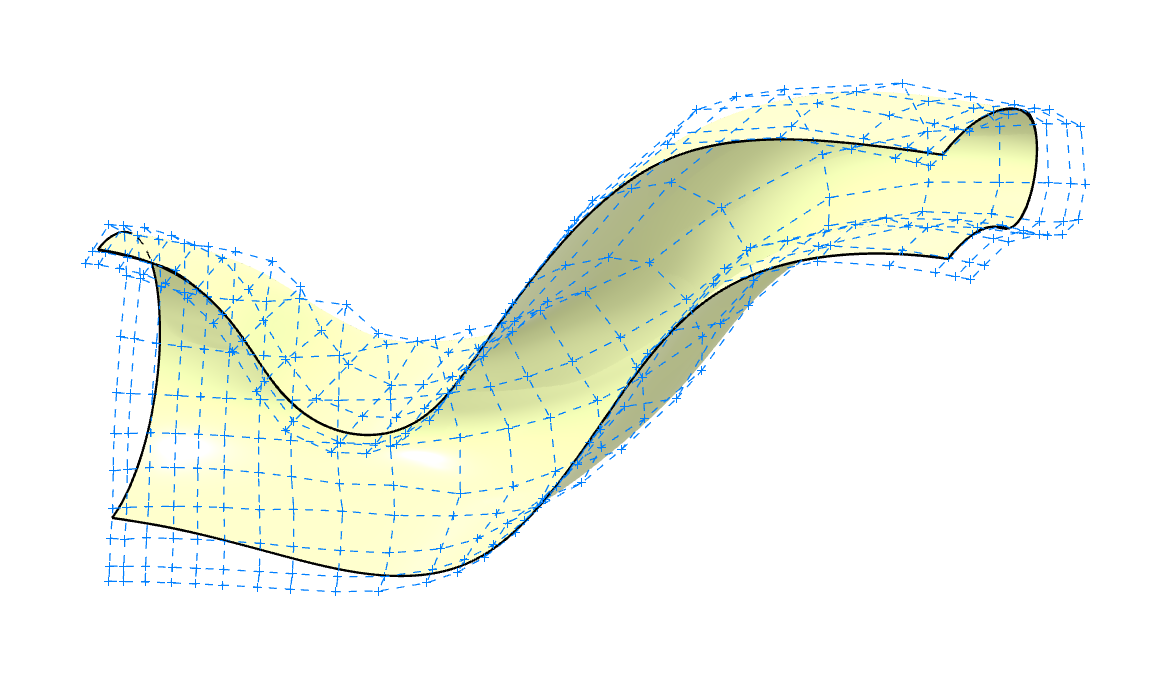}

}

\subcaption{\label{fig-meird-3dx-nurbs-param}Direct geometry
parameterisation: CATIA NURBS.}

\end{minipage}%
\begin{minipage}[t]{0.50\linewidth}

\centering{

\includegraphics[width=0.55\linewidth,height=\textheight,keepaspectratio]{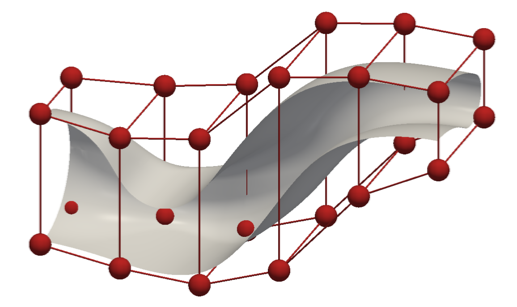}

}

\subcaption{\label{fig-meird-3dx-pygeo-param}Deformation and morphing:
FFD.}

\end{minipage}%

\caption{\label{fig-meird-geom-methods}Representative MEIRD
parameterisations for the comparator families used in the case study:
OpenVSP, CATIA B-rep, CATIA NURBS, and FFD.}

\end{figure}

The Military Engine Intake Research Duct (MEIRD) is a highly curved
three-dimensional intake surface intended to provoke the large-scale
flow distortion associated with short, strongly bent integrated inlet
systems used in unmanned air vehicles and low-observable aircraft to
reduce radar cross section (RCS) {[}33{]}. That geometry makes it a
useful constraint-editing test case because the aerodynamic path, the
packaging envelope, and the inlet and engine-face interfaces are tightly
coupled. In this study, the baseline duct was paired with a sequence of
realistic revision cases representing interface, installation, and
packaging changes. The implemented comparator families introduced in
Section~\ref{sec-related-work} can represent useful duct variation. A
reduced-order SVD comparator was not evaluated because no suitable
geometry library was available, as shown in Table~\ref{tbl-meird-setup}.
Many revised constraints nevertheless require manual changes to the
parameterisation itself, as summarised in Table~\ref{tbl-meird-burden}.
Figure~\ref{fig-meird-geom-methods} shows the representative MEIRD
parameterisations. In the modal parameterisation method (MPM), many of
the same revisions can be introduced by editing the pseudo-structure and
re-solving the associated eigenvalue problem, while the downstream
geometry mapping is held fixed.

Table~\ref{tbl-meird-burden} compares the pre-optimisation effort
required to revise the admissible duct family under seven constraint
changes. Cases 1 to 5 remain support-defined revisions, in the sense
that the revised family can still be obtained by changing supports or
releases in the pseudo-structure before the modal solve. Cases 6 and 7
remain realistic aerospace revisions, but part of the admissibility
logic then moves outside that support-definition step.

\begin{longtable}[]{@{}
  >{\raggedright\arraybackslash}p{(\linewidth - 4\tabcolsep) * \real{0.3684}}
  >{\raggedright\arraybackslash}p{(\linewidth - 4\tabcolsep) * \real{0.2632}}
  >{\raggedright\arraybackslash}p{(\linewidth - 4\tabcolsep) * \real{0.3684}}@{}}
\toprule\noalign{}
\begin{minipage}[b]{\linewidth}\raggedright
Comparator family
\end{minipage} & \begin{minipage}[b]{\linewidth}\raggedright
Implementation
\end{minipage} & \begin{minipage}[b]{\linewidth}\raggedright
Geometry input
\end{minipage} \\
\midrule\noalign{}
\endfirsthead
\toprule\noalign{}
\begin{minipage}[b]{\linewidth}\raggedright
Comparator family
\end{minipage} & \begin{minipage}[b]{\linewidth}\raggedright
Implementation
\end{minipage} & \begin{minipage}[b]{\linewidth}\raggedright
Geometry input
\end{minipage} \\
\midrule\noalign{}
\endhead
\bottomrule\noalign{}
\tabularnewline
\caption{MEIRD comparator families, implementations, and geometry
inputs.}\label{tbl-meird-setup}\tabularnewline
\endlastfoot
Compact analytic & OpenVSP {[}11{]} & OpenVSP-native paramterisations \\
Direct geometry parameterisation & CATIA B-rep & CAD sketches, splines,
and lofts \\
Direct geometry parameterisation & CATIA NURBS & B-rep converted to
NURBS \\
Deformation and morphing & FFD using pyGeo {[}34,35{]} & Triangulated
STL \\
Reduced-order & N/A & No geometry library for SVD \\
\end{longtable}

\clearpage

\begin{longtable}[]{@{}
  >{\raggedright\arraybackslash}p{(\linewidth - 8\tabcolsep) * \real{0.0588}}
  >{\raggedright\arraybackslash}p{(\linewidth - 8\tabcolsep) * \real{0.2353}}
  >{\raggedright\arraybackslash}p{(\linewidth - 8\tabcolsep) * \real{0.2353}}
  >{\raggedright\arraybackslash}p{(\linewidth - 8\tabcolsep) * \real{0.2353}}
  >{\raggedright\arraybackslash}p{(\linewidth - 8\tabcolsep) * \real{0.2353}}@{}}
\toprule\noalign{}
\begin{minipage}[b]{\linewidth}\raggedright
\#
\end{minipage} & \begin{minipage}[b]{\linewidth}\raggedright
Case
\end{minipage} & \begin{minipage}[b]{\linewidth}\raggedright
Practical driver
\end{minipage} & \begin{minipage}[b]{\linewidth}\raggedright
MPM support logic
\end{minipage} & \begin{minipage}[b]{\linewidth}\raggedright
Other Methods
\end{minipage} \\
\midrule\noalign{}
\endfirsthead
\toprule\noalign{}
\begin{minipage}[b]{\linewidth}\raggedright
\#
\end{minipage} & \begin{minipage}[b]{\linewidth}\raggedright
Case
\end{minipage} & \begin{minipage}[b]{\linewidth}\raggedright
Practical driver
\end{minipage} & \begin{minipage}[b]{\linewidth}\raggedright
MPM support logic
\end{minipage} & \begin{minipage}[b]{\linewidth}\raggedright
Other Methods
\end{minipage} \\
\midrule\noalign{}
\endhead
\bottomrule\noalign{}
\tabularnewline
\caption{Pre-optimisation effort required to revise the MEIRD admissible
family under seven representative constraint
changes.}\label{tbl-meird-burden}\tabularnewline
\endlastfoot
1 & Fixed-faces & Enforced OML and engine location & Fix inlet and
engine supports, then re-solve & Freeze or tie end controls \\
2 & Inlet-z-slide & Thrust-line or inlet-placement change & Permit the
allowed slide direction, then re-solve & Unlock controls selectively \\
3 & Added-fixture & Bracket, sensor, or duct mount & Add a local support
region, then re-solve & Lock or localise variables \\
4 & No-symmetry & Twin-engine or asymmetric installation & Remove
symmetry, then re-solve & More coupled controls \\
5 & Inlet-pin & Rotational admissibility at inlet & Add rotation at the
inlet support, then re-solve & Preserve rotation through awkward control
motion \\
6 & Stay-out-zone & Packaging or observability limit & Apply downstream
filtering or penalisation & Add exclusion constraints \\
7 & Engine-diameter-change & Changed engine selection & Rebuild
interface basis, then re-solve & Refit end controls \\
\end{longtable}

Case 1, fixed-faces, preserves the outer mould line (OML) and the
engine-face location while allowing only internal S-duct variation. This
corresponded to a common aerodynamic refinement problem in which
pressure recovery, pressure distortion, or swirl is improved without
moving either interface. The MPM representation enforced all six degrees
of freedom to zero at the inlet and engine-face supports, after which
the eigenvalue problem was re-solved to obtain the revised admissible
basis. OpenVSP, CATIA B-rep, and free-form deformation (FFD) can impose
the same practical restriction by freezing or tying the end controls.
However, the CATIA non-uniform rational B-spline (NURBS) surface is less
direct because the control points do not lie on the duct end faces, as
shown in Figure~\ref{fig-meird-3dx-nurbs-param}, so preserving those
interfaces may require surface rebuilding or reparameterisation.

Case 2, inlet-z-slide, relaxes the inlet position to allow translation
along the z axis while the remaining inlet motions stay fixed. This is
representative of installation studies driven by packaging, occlusion,
electromagnetic compatibility, or RCS requirements. In the
pseudo-structure, x and y translation and all three rotations were
retained as constraints at the inlet, while z translation was free. The
revised MPM basis therefore incorporated the permitted inlet motion
directly. OpenVSP and CATIA B-rep can represent the same revision with
an inlet-section or thrust-line z-offset variable. The CATIA NURBS and
FFD implementations are less convenient at this boundary because the
inlet is not controlled by a simple set of end-face variables and the
inlet region must be translated through grouped control motion.

Case 3, added-fixture, imposes a local installation constraint from a
bracket, sensor, or duct mount tied to a firewall or bulkhead. The
protected region is generally not aligned with the main duct sectioning,
which makes the edit awkward for parameterisations based on explicit
section controls. This is the same type of basis-design problem seen in
the wing study in Section~\ref{sec-3d-wing-basis-design}, although it is
now applied to a curved intake rather than to a shell planform. In the
MPM model, the fixture was introduced as a local support region and the
basis was updated through a new modal solve. Other methods can still
protect the same area, but they usually do so through local variable
locking, tied controls, or manual restriction of nearby geometry
handles. This is especially difficult with OpenVSP. The revised duct
family therefore remains possible, although the implementation effort
moves into local control management.

Case 4, no-symmetry, removes the symmetry assumption for a twin-engine
or otherwise asymmetric installation. Once that assumption is dropped,
the admissible duct family can shift in both y and z, and a full model
must be parameterised instead of a half model. For MPM, the change was
limited to removing the symmetry condition from the pseudo-structure and
re-computing the selected basis on the full model. The downstream
mapping from modal coordinates to geometry was unchanged. In the other
methods, the same revision commonly doubles the exposed variable count
because controls that were previously mirrored must now be managed
independently. The consequence is not loss of shape reach, but a less
compact pre-optimisation variable structure.

Case 5, inlet-pin, extends the asymmetric installation by allowing the
inlet face to rotate as well as translate. This corresponds more closely
to inlet clocking or installation-angle variation. The MPM update again
remained at the support-definition level. The inlet boundary condition
was revised to admit the required rotation, and a new basis was then
obtained from the updated eigenvalue problem. In the other methods, the
same rotational allowance must be preserved through coordinated motion
of splines, control points, or morphing boxes, which is feasible but
less direct to define and maintain.

Cases 1 to 5 therefore share a common modelling structure. The revised
design space is still defined by support conditions, releases, or local
attachments in the pseudo-structure, so the admissible family can be
updated inside the basis-definition step itself. Under those
circumstances, MPM concentrates the revision in one modelling layer:
edit the pseudo-structure, solve the eigenvalue problem again, and
retain the same downstream geometry mapping. OpenVSP, CATIA B-rep, CATIA
NURBS, and FFD can still represent useful revised ducts, but the
pre-optimisation effort more quickly shifts into explicit management of
geometric controls, protected regions, and variable couplings.

Case 6, stay-out-zone, changes the form of the constraint rather than
adding another support condition. A packaging or observability no-go
volume is an inequality constraint on the duct geometry, not a support
condition on the pseudo-structure. MPM can still generate smooth duct
families, but prohibited-volume compliance is not embedded directly in
the basis. The selected modes define admissible deformation patterns,
while clearance must be checked only after modal amplitudes have been
applied. As such, exclusion from the stay-out region must be enforced
through downstream filtering, penalisation, or rejection of infeasible
candidates. The other methods face the same restriction, although they
can often attach it more directly to local clearance checks or explicit
geometric limits possible in methods like CATIA B-rep.

Case 7, engine-diameter-change, represents a discrete downstream
interface revision caused by changed engine selection. This differs from
the earlier cases because the engine-face definition itself changes,
rather than the allowed motion within a fixed interface. MPM can still
represent the revised duct, but the pseudo-structure and selected modal
family must first be rebuilt around the new engine-face condition. The
reason is that the selected modes provide admissible shape patterns
under the chosen normalisation, not an exact prescribed interface change
in absolute length. In the other methods, the corresponding effort
typically appears as refitting or rebuilding the end controls to match
the new engine-face diameter.

The MEIRD cases show the same practical boundary seen in the earlier
examples. When the revised duct family can be expressed through
support-defined admissibility, MPM absorbs the change through
pseudo-structural edits and a new modal solve while the downstream
geometry mapping remains stable. When the governing revision becomes an
exclusion inequality or a discrete engine-face redefinition, part of the
implementation moves outside that basis-definition step. The method
remains viable, but the practical advantage becomes smaller because
clearance limits or exact interface sizes must then be enforced through
amplitude selection, downstream checks, or basis rebuilding.

\subsection{Complex Geometry: Advanced Military Airborne
Platforms}\label{sec-advanced-platforms}

\begin{figure}

\centering{

\includegraphics[width=0.8\linewidth,height=\textheight,keepaspectratio]{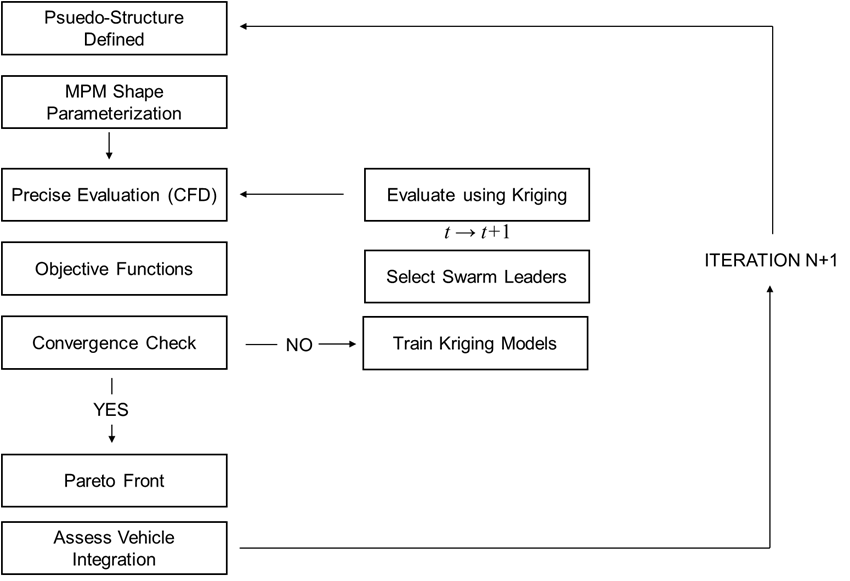}

}

\caption{\label{fig-mq-28-analysis-workflow}Workflow used in the
advanced-platform case studies, coupling modal deformation, mesh
morphing, surrogate-assisted search, and CFD.}

\end{figure}

The structural modal parameterisation method described in
Section~\ref{sec-modal-parameterisation-method} was applied to two
development problems encountered during the design of an unmanned aerial
combat vehicle (UCAV). Both cases used a consistent digital workflow
(Figure~\ref{fig-mq-28-analysis-workflow}) that coupled modal
deformation, mesh morphing, and high-fidelity CFD within a constrained
optimisation loop. The search strategy used a surrogate-assisted
Particle Swarm Optimisation (PSO) workflow consistent with the approach
of {[}36,37{]}. This nested surrogate-with-PSO pattern provided an
efficient balance between global exploration and focused local
refinement while respecting program constraints. Structural mode shapes
were computed with ANSYS Mechanical™ and aerodynamic analyses were
performed with ANSYS Fluent™. Projection of structural modal motion to
the CFD OML used the Modal Projection and Force Reconstruction (MPR)
formulation of {[}31{]}.

\subsubsection{Case 1: Bifurcated Serpentine
Intake}\label{case-1-bifurcated-serpentine-intake}

The first case addressed a trade study for a bifurcated serpentine
intake geometry at the conceptual design phase. The intake had to
satisfy both aerodynamic performance requirements, namely preserving
mass flow and minimising fan-face distortion, and stringent geometric
exclusion zones imposed by internal payload and signature requirements.
The objectives were defined according to engine-manufacturer
specifications for ideal installation (zero losses).

The baseline concept was evaluated at 12 mission-relevant flight
conditions to establish a reference point. Given operational
requirements, the intake performance was considered as a multi-point,
multi-objective problem. To limit the number of computationally
expensive high-fidelity CFD simulations, only three design reference
conditions (shown in Table~\ref{tbl-mq-28-design-ref-conds}) were
selected for the optimisation loop. These conditions were identified to
capture relevant mission trade-offs and balance pressure recovery and
distortion metrics for a critical static-thrust case (ID 1B),
high-angle-of-attack low-altitude airfield departure (ID 7), and
low-angle-of-attack high-altitude performance (ID 11).

\begin{longtable}[]{@{}lcrc@{}}
\toprule\noalign{}
ID & Mach & Altitude & Angle of Attack \\
\midrule\noalign{}
\endfirsthead
\toprule\noalign{}
ID & Mach & Altitude & Angle of Attack \\
\midrule\noalign{}
\endhead
\bottomrule\noalign{}
\tabularnewline
\caption{MQ-28 intake design reference conditions retained for the
optimisation loop.}\label{tbl-mq-28-design-ref-conds}\tabularnewline
\endlastfoot
1B & 0 & Sea Level & 0 \\
7 & Low subsonic & Sea Level & Near-stall \\
11 & Transonic & Ceiling & Cruise \\
\end{longtable}

Aerodynamic performance was evaluated in CFD using steady RANS with the
\(k-\omega\) SST turbulence model on a volume mesh of approximately 20
million cells. Mesh refinement targeted the inlet boundary layer so that
\(y^+ \approx 1\) on intake walls. The pseudo-structure for the MPM was
developed directly from the CFD mesh. The intake lip was assumed to be a
fully fixed connection at the entrance to the duct. The Aerodynamic
Interface Plane (AIP) was treated as a sliding support along the x-axis.

The first seven low-order modes, comprising three bending, two
torsional, and two breathing modes, were evaluated. Their amplitudes
formed the compact design vector \(a \in \mathbb{R}^7\). Payload and
low-observable constraints were imposed only through limiting the
design-variable range. In the final solution, these constraints
restricted the AIP x-location design-variable range such that a fully
fixed connection in the pseudo-structure would likely have produced a
better initial representation.

Modal shape deformations (\(\Phi\)) were projected to the CFD OML using
the MPR transfer operator from {[}31{]} with nearest-neighbour
connectivity (\(L = 4\)), producing a sparse transfer matrix (\(T\)).
The deformation Jacobian (\(J_\text{def} = T \Phi\)) was exported and
archived for adjoint coupling and gradient verification. This enabled
efficient mapping of Fluent OML sensitivities back to modal amplitudes
when adjoint data became available.

The optimisation was bound by constraints to maintain \(\geq\) 95\%
operating efficiency at conditions ID 7 and ID 11
(Table~\ref{tbl-mq-28-design-ref-conds}), with two design objectives
targeting full unity (ideal recovery) at the AIP and minimising the
forward loading of the duct (representing system head loss). Static
pressure was imposed as a CFD boundary condition, so mass flow through
the intake remained an implicit result; this was considered necessary to
promote convergence of both the CFD solution and the optimisation
process. The design search used a surrogate-assisted multi-objective PSO
to find promising trade regions on Kriging surrogates, then selected
candidates were subjected to batched high-fidelity CFD evaluation per
Figure~\ref{fig-mq-28-analysis-workflow}. A local gradient-based
optimiser was used for local refinement where adjoint sensitivities were
computed; otherwise, central finite differences
(\(\varepsilon = 1 \times 10^{-3}\)) were used.

\begin{figure}

\centering{

\includegraphics[width=0.75\linewidth,height=\textheight,keepaspectratio]{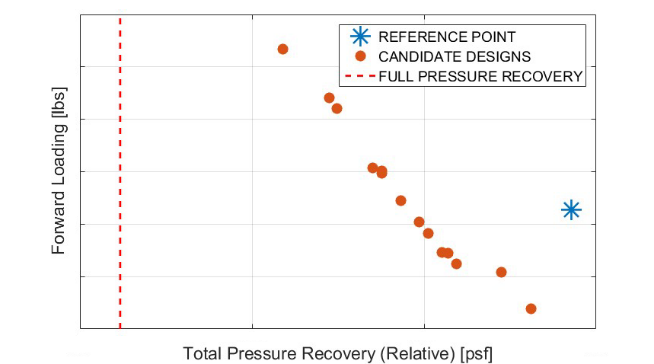}

}

\caption{\label{fig-mq-28-intake-pareto-front}Final intake Pareto front
for the MQ-28 case. The circled region marks the selected design
neighbourhood; axis values are redacted.}

\end{figure}

\begin{figure}

\centering{

\includegraphics[width=0.6\linewidth,height=\textheight,keepaspectratio]{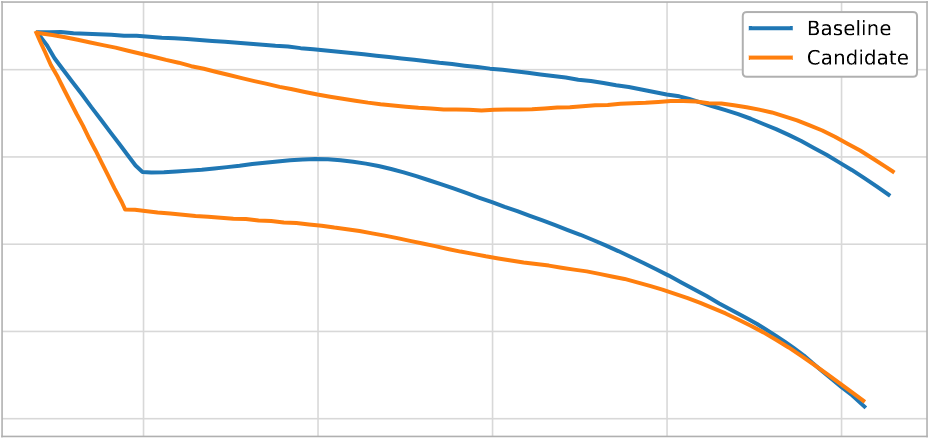}

}

\caption{\label{fig-mq-28-intake-mid-plane}Selected MQ-28 intake
mid-plane compared with the baseline; dimensions are redacted.}

\end{figure}

\begin{figure}

\centering{

\includegraphics[width=0.8\linewidth,height=\textheight,keepaspectratio]{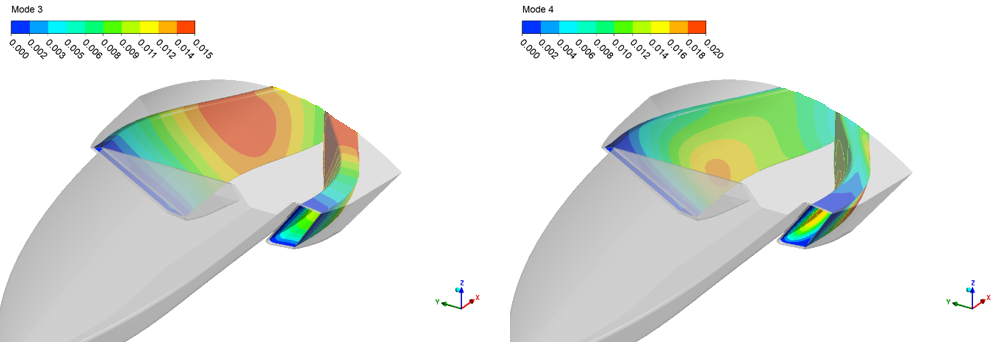}

}

\caption{\label{fig-mq-28-intake-modes}Selected MQ-28 intake control
modes. Mode 3 primarily opens and straightens the flow path, while Mode
4 provides torsional control near the high-curvature region.}

\end{figure}

Figure~\ref{fig-mq-28-intake-pareto-front} presents the final Pareto
front, and Figure~\ref{fig-mq-28-intake-mid-plane} shows the
corresponding mid-plane geometry. Convergence required three MPM update
iterations to impose additional constraints and progressively reduce
free-form deformation {[}14{]}. The selected design
(Figure~\ref{fig-mq-28-intake-modes}) showed an improvement in total
relative pressure recovery of about 12\% while respecting exclusion
zones. System head loss (forward loading) was also reduced by about 7\%.

\subsubsection{Case 2: Nose Pitch (Trim / Longitudinal Stability
Reconciliation)}\label{case-2-nose-pitch-trim-longitudinal-stability-reconciliation}

The second case study described an application used when mature layout
releases produced an adverse shift in the centre of gravity that created
a longitudinal stability risk. Large configuration changes were
impractical due to schedule and structural maturity. The nose geometry
was selected for MPM optimisation as a local, low-impact control region
to restore trim and stability margins while minimising downstream design
rework.

The baseline nose OML was evaluated at three conditions across the
flight envelope, with a coarse CFD model (approximately 12 million
cells) enabled by pressure-dominant design objectives and accelerated
numerical evaluation in support of program schedules. A steady RANS
solution provided the aerodynamic reference state for the baseline
design.

The pseudo-structure for MPM was assembled directly on the CFD surface:
shell/beam elements were attached to CFD surface nodes so that modal
DOFs mapped cleanly to aerodynamic nodes used in projection. Generic
steel material properties were assigned, and a fixed boundary (all DOFs
= 0) was applied at the rear structural frame attachment to represent
the primary support. The structure was constrained against torsional
rotation to prevent asymmetry. The upper skin was stiffened to reduce
adverse movement. To tune modal solutions and suppress content that
conflicted with integration to the existing belly substructure, a rigid
reference point was placed slightly ahead of the aft frame on the lower
skin. The tuned finite-element pseudo-structure was used to solve the
eigenvalue problem and extract candidate normal modes.

Eight eigen solutions were computed and inspected. Modes 1 and 2
exhibited global pillowing and bending and were selected for low-order
aerodynamic control. Modes 3-8 contained higher spatial-frequency local
flexures that would increase design dimensionality and risk generating
small geometric features difficult to mesh robustly. Generally, these
were excluded from the basis, although the localised deformation in mode
5 provided a useful control effect on forebody pressure and trim. This
localised feature was extracted through node indexation and added to the
selected basis. Mode selection included verification of orthogonality
and visualisation of surface-normal displacement.

The final shape parameterisation basis for optimisation comprised modes
1, 2, and the localised component derived from mode 5. All three basis
directions were normalised and implemented directly as optimiser
coordinates. Practical parameter bounds were imposed to ensure
manufacturability and solver robustness. The MPR transfer operator
{[}31{]} was validated via small perturbation tests.

Aerodynamic evaluations used Fluent steady RANS with the
\(k\)-\(\omega\) SST turbulence model. The optimisation objective
combined a trim metric (minimise required elevator deflection) and an
aerodynamic penalty (drag increase), yielding a small multi-objective
problem. The optimiser enforced hard geometric limits to maintain
internal volume, and a manufacturability filter set minimum radii on
skin edges.

\begin{figure}

\centering{

\includegraphics[width=0.8\linewidth,height=\textheight,keepaspectratio]{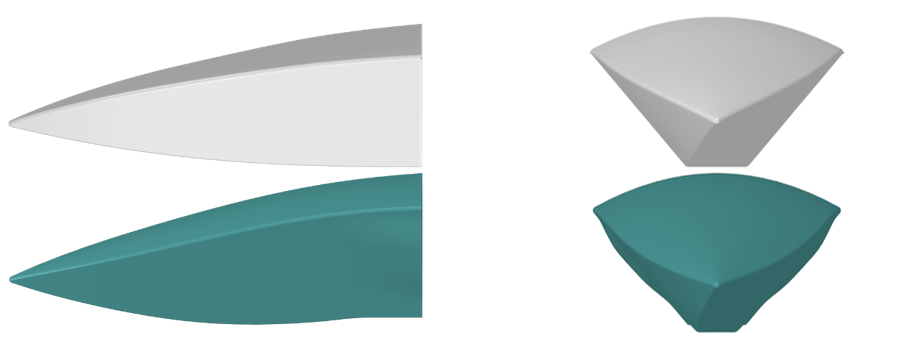}

}

\caption{\label{fig-mq-28-nose-designs}Baseline (top) and optimised
(bottom) nose geometries for the trim-reconciliation case.}

\end{figure}%

The nose modification, shown in Figure~\ref{fig-mq-28-nose-designs},
reduced required elevator trim by an average of about 0.5 degrees across
the analysed conditions while producing a net drag change of less than
1\%. Convergence was achieved in a single MPM iteration. Given that
modal amplitudes were small, mesh quality remained acceptable throughout
the process, with smoothing and dynamic mesh updates only.

\section{Discussion}\label{sec-discussion}

Structural modal parameterisation is most compelling when treated as a
parameterisation choice for constrained aerodynamic design rather than
as a universal optimisation strategy. Its practical value comes from the
selected deformation family defined before optimisation through
supports, coupling assumptions, and basis-shaping property edits. This
gives concise coordinates for search, while lower-order selected modes
often provide recognisable global shape change and higher-order selected
modes provide local refinement. The selected basis also produces smooth
deformation fields rather than pointwise geometric edits.

The method is particularly well suited to tightly constrained
optimisation problems. Once packaging constraints, manufacturability,
systems integration, or mission requirements significantly restrict the
allowable design space, exposing a large number of arbitrary geometry
controls can become expensive and difficult to manage. The beam
comparison showed that alignment with the desired admissible family
could matter more than unrestricted shape reach. The thin-wing study
showed that stiffness edits, attached masses, and protected hinge
constraints could shift the selected modal family toward the required
aerodynamic response before optimisation began. The MEIRD cases showed
that many support-defined admissibility revisions could be absorbed by
editing the pseudo-structure and re-solving the eigenvalue problem while
the downstream geometry mapping was unaffected. The advanced-platform
applications showed that the same approach remained viable on more
complex aerospace geometry when the required motion was smooth, coupled,
and strongly constrained by program requirements.

The role of optional QR conditioning should be interpreted in the same
practical way. Its purpose is numerical conditioning rather than basis
definition. It may improve numerical efficiency by reducing coordinate
coupling for a chosen optimisation algorithm. However, it does not alter
the admissible deformation family. It only changes the coordinates used
to describe the same selected subspace.

The limitations are also reasonably clear. Smooth modal deformation is
less suitable for geometries characterised by sharp corners, hard planar
breaks, or intentional discontinuities. In such cases, many selected
modes may be required before the required feature is represented
adequately, which reduces the conciseness advantage. The same limitation
applies to large planform edits such as a new cranked wing or a large
sweep break, where a smooth modal basis is better suited to refining an
existing layout than to creating a new one concisely. The method is also
sensitive to pseudo-structural modelling choices, especially boundary
conditions, selected mode count, amplitude bounds, and the basis-shaping
property edits. These are not minor implementation details. These
choices are part of the parameterisation definition itself.

Scope therefore remains important. This work does not claim that
structural modal parameterisation is the best geometry basis for every
aerodynamic optimisation problem. The beam comparison showed that
trained reduction can outperform a pseudo-structure-derived basis once
the desired family has already been generated, because the trained basis
is then specialised to that family. The thin-wing study showed that
protected admissibility can require a deliberate trade between
hinge-line coherence and whole-field flexibility. Instead, the method
offers a competitive, engineer-usable balance of conciseness, constraint
awareness, and usable flexibility for tightly constrained problems.
Within that scope, the method represents a practical and defensible
choice.

\section{Conclusion}\label{sec-conclusion}

A structural modal parameterisation method has been presented for
aerodynamic shape optimisation problems in which the allowable geometry
change is small, smooth, and strongly constrained. The method constructs
a tunable pseudo-structure, solves its eigenvalue problem, and uses
selected modes as geometry coordinates. Boundary conditions, coupling
assumptions, stiffness, density, thickness, and added masses are treated
as basis-design choices fixed before optimisation, so a substantial part
of the admissible deformation family is defined before the search
begins.

The reported studies demonstrate that this basis-design step strongly
influences the quality of the selected coordinates. The beam comparison
showed that a pseudo-structure-derived basis can recover useful coupled
motion without a precomputed shape library. The thin-wing study showed
that stiffness tuning, attached masses, and hinge constraints shifted
the selected modal family toward the required aerodynamic response while
reducing the number of selected modes. The MEIRD and advanced-platform
applications showed that the same approach extends to complex aerospace
geometry, where many support-defined revisions can be absorbed through
pseudo-structural edits and a new modal solve while the downstream
geometry transfer remains fixed.

For constrained aerodynamic design, another practical value is that the
admissible deformation family can be shaped and inspected before
optimisation. The resulting coordinates remain solver-agnostic, concise,
smooth, and engineer-interpretable, which is useful when payload
requirements, manufacturability limits, and low-observability
considerations have already reduced the available design space.

The method remains limited when the required change contains sharp local
features, inequality volume constraints, or discrete interface
redefinitions, because smooth modal bases then become less concise and
part of the admissibility logic must be enforced outside the
basis-definition step. Performance also depends directly on
pseudo-structural modelling choices, selected mode count, and amplitude
bounds. Future work should examine automated basis tuning, stronger
coupling to adjoint-based optimisation, and hybrid global-local
parameterisations for industrial design problems.

\section*{Copyright Statement}\label{sec-copyright}
\addcontentsline{toc}{section}{Copyright Statement}

The authors confirm that they, and/or their company or organisation,
hold copyright on all original material included in this paper. The
authors also confirm that they have obtained permission from the
copyright holder of any third-party material included in this paper to
publish it as part of their paper. The authors confirm that they give
permission, or have obtained permission from the copyright holder of
this paper, for the publication and distribution of this paper as part
of the ICAS proceedings or as individual off-prints from the
proceedings.

\section*{References}\label{references}
\addcontentsline{toc}{section}{References}

\protect\phantomsection\label{refs}

\end{document}